\documentclass[12pt]{amsart}
\usepackage[]{amsmath, amsthm, amsfonts, verbatim, amssymb, mathrsfs, mathtools, booktabs, multirow}

\usepackage[all]{xy}

\usepackage{hyperref}
\usepackage{pifont}

\usepackage[normalem]{ulem} 

\usepackage{enumerate}

\usepackage{graphicx}
\usepackage{epstopdf}
\DeclareMathAlphabet{\mathpzc}{OT1}{pzc}{m}{it}
\usepackage{tikz,tikz-cd, color}
\usepackage[makeroom]{cancel}
\usetikzlibrary{arrows}
\usepackage{microtype}
\usepackage{mathdots}

\usepackage{float}
\usepackage{caption, subcaption}

\usepackage{etoolbox}
\ifdef{\crop}{%
\usepackage[includeheadfoot,twoside=False,paperwidth=448pt,paperheight=587pt,rmargin=15pt,lmargin=15pt,tmargin=15pt,bmargin=15pt]{geometry}%
}{%
\setlength{\topmargin}{20mm}
\addtolength{\topmargin}{-1in}
\setlength{\oddsidemargin}{27mm}
\addtolength{\oddsidemargin}{-1in}
\setlength{\evensidemargin}{27mm}
\addtolength{\evensidemargin}{-1in}
\setlength{\textwidth}{156mm}
\setlength{\textheight}{242mm}
}%

\tikzset{
  symbol/.style={
    draw=none,
    every to/.append style={
      edge node={node [sloped, allow upside down, auto=false]{$#1$}}}
  }
}

\newtheorem{theorem}{Theorem}[section]
\newtheorem{lemma}[theorem]{Lemma}
\newtheorem{proposition}[theorem]{Proposition}

\newtheorem{conjecture}[theorem]{Conjecture}
\newtheorem{claim}[theorem]{Claim}
\newtheorem{definition}[theorem]{Definition}

\theoremstyle{remark}
\newtheorem{remark}[theorem]{Remark}

\theoremstyle{definition}
\newtheorem{notation}[theorem]{Notation}

\newcommand{\bC}{\mathbb{C}}

\newcommand{\bQ}{\mathbb{Q}}
\newcommand{\bP}{\mathbb{P}}
\newcommand{\bR}{\mathbb{R}}
\newcommand{\bZ}{\mathbb{Z}}

\newcommand{\cC}{\mathcal{C}}

\newcommand{\cE}{\mathcal{E}}
\newcommand{\cF}{\mathcal{F}}

\newcommand{\cH}{\mathcal{H}}

\newcommand{\cO}{\mathcal{O}}

\newcommand{\cQ}{\mathcal{Q}}
\newcommand{\cS}{\mathcal{S}}
\newcommand{\cT}{\mathcal{T}}
\newcommand{\cU}{\mathcal{U}}
\newcommand{\cV}{\mathcal{V}}

\newcommand{\tbP}{\widetilde{\mathbb{P}}}

\newcommand{\tcS}{\widetilde{\mathcal{S}}}

\newcommand{\tX}{\widetilde{X}}
\newcommand{\tY}{\widetilde{Y}}

\newcommand{\id}{\mathrm{id}}

\newcommand{\cMov}{\overline{\mathrm{Mov}}}
\newcommand{\Bir}{{\rm Bir}}
\newcommand{\Aut}{{\rm Aut}}

\newcommand{\Bl}{{\rm Bl}}

\newcommand{\coker}{{\rm coker}\,}

\newcommand{\Eff}{\mathrm{Eff}}

\newcommand{\rBs}{\mathrm{Bs}}

\newcommand{\Exc}{{\rm Exc}}

\newcommand{\Gr}{{\rm  Gr}}

\newcommand{\Hom}{{\rm Hom}}

\newcommand{\sHom}{\mathscr{H}om}

\newcommand{\Pic}{{\rm Pic}\,}

\DeclareMathOperator{\Nef}{Nef}
\DeclareMathOperator{\Sing}{Sing}
\DeclareMathOperator{\rank}{rank}

\def\<{\langle}
\def\>{\rangle}

\def\ra{\rightarrow}
\def\dra{\dashrightarrow}

\address[]{Department of Mathematics, 
Faculty of 
Environmental, Life, Natural Science and Technology,
Okayama University,
Okayama, Japan
}
\email{ito-atsushi@okayama-u.ac.jp}

\address[]{Department of Mathematics, National Cheng Kung University, Tainan, Taiwan
}
\email{cjlai72@mail.ncku.edu.tw}

\address[]{Institute of Mathematics, Academia Sinica, Taipei, Taiwan
}
\email{sswangtw@gate.sinica.edu.tw}


\begin{document}
\title[The movable cone of CY3 in ruled Fano manifolds]{The movable cone of Calabi--Yau threefolds in ruled Fano manifolds}

\author{Atsushi\ Ito,\ Ching-Jui\ Lai, \ Sz-Sheng Wang}

\begin{abstract} 
We describe explicitly the chamber structure of the movable cone for a general complete intersection Calabi--Yau threefold in a non-split $(n + 4)$-dimensional $\mathbb{P}^{n}$-ruled Fano manifold of index $n + 1$ and Picard number two. Moreover, all birational minimal models of such Calabi--Yau threefolds are found whose number is finite.
\end{abstract}

\maketitle


\section{Introduction}\label{intro_sec}

In the classical Mirror Symmetry, the mirror map conjecturally identifies a neighborhood of a large volume limit point in the K\"ahler moduli of a smooth Calabi--Yau threefold $X$ with that of a special boundary point in the complex moduli of a mirror, which is characterized by unipotent monodromy. If the nef cone $\Nef (X)$ is a finite rational polyhedral cone, there is a partial compactification with a finite number of large volume limit points of the complexified K\"ahler cone of $X$. In general, a cone conjecture proposed by Morrison and Kawamata \cite{Morrison93,Morrison96,Kawamata97} predicts that the nef cone has only finitely many $\Aut(X)$-orbits of edges.

In the present article, by a smooth Calabi--Yau threefold, we mean a smooth projective threefold $X$ with $K_X \sim 0$ and $H^1 (\cO_X) = 0$. We focus on the birational version of the cone conjecture for such $X$.  Recall that a divisor $D$ is \emph{movable} if for some positive number $m$ the linear system $|m D|$ has no fixed components 
and the \emph{movable cone} $\cMov (X)$ is the closure of the convex hull of movable divisor classes. The movable cone conjecture is the following:

\begin{conjecture}[\cite{Morrison96,Kawamata97}]\label{coneconj}
There is a finite rational polyhedral cone which is a fundamental domain for the action of the birational automorphism group $\Bir (X)$ on the movable effective cone $\cMov (X) \cap \Eff (X)$.
\end{conjecture}

For a survey of this widely open conjecture, we refer the reader to \cite{Totaro10, LOP18}.

A \emph{normalized $\bP^n$-ruled Fano manifold $P$ over $M$} is the projective bundle associated to an ample bundle $\cF$ of rank $n + 1$ with $c_1 (\cF) = c_1 (T_M)$. Such a pair $(M, \cF)$ is called a \emph{Mukai pair} in \cite{Kan19A}. The classification of Mukai pairs with $\rank(\cF)\geqslant\dim M-2$ has been completed recently in \cite{Kan19}, and in this series of work we mainly consider the case when $\dim M=4$ and Picard number $\rho(M)=1$. When $\cF$ splits as a direct sum of line bundles, we proved Conjecture \ref{coneconj} in \cite{LW2022} for a general complete intersection Calabi--Yau threefold in such $P$ except when $M$ is del Pezzo of degree one. 
In this paper, we treat the non-splitting $\cF$ as follows,
see Theorem \ref{mainV4}, \ref{mainV5} and \ref{mainG24} for the details.

\begin{theorem} \label{mainthm_intro}
Let $P = \bP (\cF)$ be a normalized $\bP^n$-ruled Fano manifold over a smooth projective fourfold $M$ with $\rho (M) = 1$. Assume that $\cF$ is not isomorphic to a direct sum of line bundles nor $T_{\bP^4}$. Then a complete intersection $X_{\cF}$ of $n + 1$ general hypersurfaces in $|\cO_P (1)|$ is a smooth Calabi--Yau threefold of Picard number two. Moreover, all the minimal models of $X_\cF$ are constructed and the movable cone $\cMov (X_{\cF})$ is a  
rational polyhedral cone.
\end{theorem}

Note that the case $\cF = T_{\bP^4}$ is non-splitting but it can be reduced to the case $\cO_{\bP^4} (1)^{\oplus 5}$ (see \cite[\S6.1]{LW2022}). Therefore, Theorem \ref{mainthm_intro} does imply that Conjecture \ref{coneconj} holds for a general complete intersection Calabi--Yau threefold in a normalized $\bP^n$-ruled Fano manifold of Picard number two associated to a non-splitting bundle over a smooth Fano fourfold.

The proof of Theorem \ref{mainthm_intro} is based on the construction of the determinantal contraction in \cite{SSW18} and then modeled on the proof for splitting $\cF$ cases in \cite{LW2022}, but requires many new ideas. We briefly discuss the technical issues below.

For $\cF$ being splitting, we can find a special surface in $X_\cF$ by applying the geometric construction of Eagon--Northcott complexes (see \cite[Proposition 4.4]{LW2022}). When $\Bir (X_\cF)$ is finite, studying the geometry of such surfaces is  sufficient for finding all the other minimal models of $X_\cF$. However, this construction does not work when $\cF$ is non-splitting.

When $\cF$ is not isomorphic to $T_{\bP^4}$ and non-splitting, there are three cases (see Theorem \ref{clFano}), which are bundles on the Grassmannian $\Gr (2, 4)$, and del Pezzo fourfolds $V_4$ and $V_5$ of degrees 4 and 5 respectively. Except for the case $V_4$, the other two cases happen to have very different technical difficulties from the splitting cases:

In each case, we have the associated small contraction from $X_{\cF}$ to a determinantal variety in $M$ and its flop $X_{\cF} \dashrightarrow X_{\cE} $ (see \eqref{flopdiag}).
Since there is a natural fibration $X_{\cE} \to \bP^n$ with $n=1$ or $2$,
we obtain one edge of $\cMov (X_{\cF})$.
Thus what we need to do first is to find another contraction $X_{\cF} \to Y_{\cF}$ and its flop.

\begin{enumerate}[(i)]
    \item The constructions of $X_{\cF} \to Y_{\cF}$ and its flop in the case $M = V_4$ are similar to those in the case $(\Gr(2,4),\cO(1)^{\oplus 4})$ (see \cite[\S 6.2]{LW2022}). 
    By construction, we have a generically $2 : 1$ morphism $X_{\cF} \to \bP^3$. 
    We can check that its stein factorization $X_\cF \to Y_\cF$ is a small contraction. 
    Since there is an involution $ X_{\cF} \dashrightarrow X_{\cF}$ over $\bP^3$, which is the flop of $X_{\cF} \to Y_{\cF}$, we obtain another edge of $\cMov (X_{\cF})$. Thus we can determine $\cMov (X_{\cF})$ and its decomposition by nef cones of all minimal birational models (Theorem \ref{mainV4}). 
    
    \item\label{V5_introsec} For the case $M = V_5$, we give an alternative description of $X_{\cF}$ as a suitable zero locus in $\bP_{\bP^4}(\wedge^3T_{ \bP^4})$. Moreover, by finding out an analog of $(\cE, \cF)$, we can apply the determinantal construction of \cite{SSW18} to this description to obtain a small contraction $X_\cF \to Y_\cF$ and its flop $X_{\cF} \dashrightarrow X_\cF^+$. 
    Finally, using the method in \cite{Ito14} of constructing Mori dream spaces with Picard number two, we can construct a small contraction of $X_\cF^+$, its flop $X_\cF^+ \dashrightarrow X_\cF^{++}$ and a fibration $X_\cF^{++} \to \bP^1$.
    Thus we can determine $\cMov (X_{\cF})$ and its decomposition (Theorem \ref{mainV5}).
    
    \item For the case $M = \Gr (2, 4)$, we give an alternative description of $X_{\cF}$ as in \eqref{V5_introsec}, which induces a small contraction $X_\cF \to Y_\cF \subseteq \bP^4$. 
    However, this description involves a morphism $\tau$ from a vector bundle 
    to a reflexive sheaf\footnote{It is called a B\u anic\u a sheaf, see \cite[Definition 1.5]{Kan19} and references therein.}, 
    which is not locally free at a unique point $p \in \bP^4$  (see \eqref{G24XF+}). 
    Hence we ``resolve the singularity'' of $\tau$, that is, consider a relevant morphism $\varrho$ between vector bundles  
    on the blow-up $\tbP \to \bP^4$ at $p$ (see \eqref{G24rho}). By the determinantal construction, the morphism $\varrho$ induces a small contraction $\tX_\cF\ra\tY_\cF \subseteq \tbP$, which descends to the flop $X_{\cF}^{+} \to Y_\cF$ of $X_\cF \to Y_\cF$. Since there is a fibration $X_{\cF}^{+} \to \bP^2$,
    we obtain another edge of $\cMov (X_{\cF})$ (Theorem \ref{mainG24}).
\end{enumerate}

We remark that Conjecture \ref{coneconj} has been verified for several special cases, see \cite{Borcea91, Kawamata97, Fryers01, LP13, Og-CY2, CO2015, BorisovNuer2016,yanez2022,WL2021remarks,LW2022} and references therein. For the case of Calabi--Yau manifold $X$ with Picard number two, Lazi\'c and Peternell \cite[Theorem 1.4]{LP13} proved that there is a polyhedral cone $\Pi$ which is a fundamental domain for the action of $\Bir (X)$ on the movable effective cone of $X$, and that the cone $\Pi$ is rational under the additional assumption that $\Bir (X)$ is infinite, i.e., Conjecture \ref{coneconj} holds for such $X$. However, Conjecture \ref{coneconj} is still open when $\Bir (X)$ is finite. The fundamental difficulty is whether $\cMov (X)$ is rational. Our main result provides new positive evidence in this case.

The paper is organized as follows: In \S\ref{pre_sec} we recall some foundational material concerning the classification of Fano bundles and the determinantal construction. We also give Table \ref{intAll} and \ref{topAll} of intersection numbers and Hodge numbers on our Calabi--Yau threefolds respectively. We devote the remaining sections \S\ref{V4_sec}-\S\ref{G24_sec} to the proof of our main result, Theorem \ref{mainthm_intro}.

\begin{notation}
Throughout this paper, we work over the complex field $\bC$. For a coherent sheaf $\cF$ on a variety $M$ and an integer $r \geqslant 1$, we write $\Gr_M (\cF, r)$ for the Grassmannian of $r$-dimensional quotients of $\cF$. For the case $r = 1$, we denote by $\bP_M (\cF)$ the projective bundle with the tautological line bundle $\cO_{\cF}(1) \coloneqq \cO_{\bP(\cF)} (1)$. When $\cF$ is a vector bundle of rank $f$, we set $\Gr_M (r, \cF) \coloneqq \Gr_M (\cF, f - r)$. Moreover, we will use the notation $\Gr_M (r ,f)$ if $\cF$ is trivial and omit the subscript ``$M$'' when no confusion can arise.
\end{notation}

\medskip

\emph{Acknowledgements.} 
A.\ Ito was supported by JSPS KAKENHI Grant Number 17K14162, 21K03201. C.-J.\ Lai is supported by the National Center for Theoretical Sciences in Taiwan. Part of this work was done when he visited A.\ Ito at Nagoya University and he thanks the institute for the warm hospitality. C.-J.\ Lai and S.-S.\ Wang are supported by the National Sciences and Technology Council in Taiwan under grant number 109-2115-M-006-013-MY2 and 111-2115-M-001-003-MY3 respectively. S.-S.\ Wang thanks the Institute of Mathematics at Academia Sinica for providing support and a stimulating environment.



\section{Preliminaries}\label{pre_sec}

\subsection{Non-split Fano bundles} \label{fb_subsec}

Let us briefly recall the classification of Fano bundles, i.e., vector bundles whose projectivization are Fano manifolds. For simplicity, we only list such bundles which will be needed in this paper.

\begin{theorem}[\cite{PSW92,Occ01,NO07, Kan19}] \label{clFano}
Let $\cF$ be an ample bundle on a smooth projective fourfold $M$ with $c_1 (\cF) = c_1 (T_M)$ and $\rho (M) = 1$. If $\cF$ is not isomorphic to a direct sum of line bundles nor to $T_{\bP^4}$, then the pair $(M , \cF)$ is isomorphic to one of the following:
\begin{center}
    $(V_4, p^{\ast} \cS^{\vee} (1))$, \quad $(V_5, \cS^{\vee}_{V_5}(1))$ \quad \text{or} \quad $(\Gr(2, 4), \cS (2) \oplus \cO (1))$.
\end{center}
Here we use the following symbols:
\begin{itemize}
    \item $\cS$ is the universal subbundle on $\Gr (2, 4)$.
    \item $V_4$ is a quartic del Pezzo fourfold obtained as a double cover $V_4 \xrightarrow{p} \Gr(2, 4) \cong Q_4$ of the hyperquadric of dimension $4$ branched along a smooth divisor $B \in |\cO_{Q_4} (2)|$.
    \item $V_5$ is a general $2$-codimensional linear section of the Grassmannian $\Gr(2, 5)$ embedded into the projective space $\bP^9$ via the Pl\"ucker embedding.
    \item $\cS_{V_5}$ is the restriction of the universal subbundle on $\Gr (2, 5)$ to $V_5$.
\end{itemize}
\end{theorem}

The Fano bundle appearing above defines a normalized Fano manifold ruled by projective spaces. Let us recall the definition of such manifolds from \cite[Definition 3.2]{NO07}.

\begin{definition} \label{ruledFano}
Let $\cF$ be a vector bundle of rank $r \geqslant 2$ on a projective manifold $M$. We will call $\bP (\cF)$ a normalized $\bP^{r - 1}$-ruled Fano manifold if $\cF$ is an ample bundle with $c_1 (\cF) = c_1 (T_M)$.
\end{definition}

Note that such Fano manifold $\bP (\cF)$ has index $r$, i.e., $\cO (K_{\bP (\cF)}) \cong  \cO_\cF (- r)$. Conversely, each $\bP^{r - 1}$-ruled Fano manifold $\bP (\cF)$ of index $r$ is isomorphic to a normalized one (see \cite[Proposition 3.3]{NO07}).

\subsection{Determinantal varieties} \label{detctrs_subsec}

This subsection reviews briefly some results of degeneracy loci, which will be needed in the later sections. For more details we refer the reader to \cite[\S3 and \S4]{LW2022}.

Fix a morphism $\sigma \colon \cE^{\vee} \to \cF$ of vector bundles on a variety $M$ of rank $e$ and $f$ respectively. For each $k \leqslant \min \{e, f\}$, the $k$th \emph{degeneracy locus} of $\sigma$ is
\[
    D_k (\sigma) = \{ x \in M \mid \rank (\sigma (x)) \leqslant k\},
\]
with the convention $D_{- 1} (\sigma) = \varnothing$. It can be described as the zero locus of the global section $\wedge^k \sigma$ of the bundle $\wedge^{k + 1} \cE \otimes \wedge^{k + 1} \cF$. Note that the \emph{expected codimension} of $D_k (\sigma)$ in $M$ is $(e - k) (f - k)$, though $D_k(\sigma)$ can be empty or have strictly smaller codimension.
However, it is known that $D_k (\sigma) = \varnothing$ or the codimension of $D_k (\sigma)$ in $M$ coincides with the expected one
if $\cE \otimes \cF$ is globally generated and $\sigma$ is general (see \cite{Banica91} for example).

The following lemma will be needed in Section \ref{G24_sec}.

\begin{lemma}\label{eq_as zero locus}
For an integer $r \geqslant 1$, consider the Grassmannian $\pi \colon \Gr (\cF, r) \ra M$ as a Quot scheme. Let $\pi^* \cF \ra \cQ$ be the universal quotient bundle of rank $r$ on $\Gr(\cF, r)$ and $\cC = \coker \sigma$. Then the closed subscheme
\begin{align}
        \Gr (\cC, r) \subseteq \Gr (\cF,r) 
\end{align}
is the zero locus of the composition $\pi^* \cE^{\vee} \xrightarrow{\pi^\ast \sigma} \pi^* \cF \ra \cQ$ in $\Gr (\cF,r)$.
\end{lemma}

\begin{proof}
Recall that a morphism $g \colon T \ra  \Gr (\cF ,r )$ corresponds to a quotient bundle $ h^* \cF \ra Q$ of rank $r$ over $h \colon T \ra M$. In this correspondence, $ h^* \cF \ra Q $ coincides with the pullback of $\pi^* \cF \ra \cQ $ by $g$.
Then $g$ factors through $  \Gr (\cC,r) \hookrightarrow  \Gr (\cF,r) $ if and only if $  h^* \cF \ra Q$ factors through $h^* \cF \ra  h^* \cC $ if and only if the composite 
$h^* \cE^{\vee} \ra h^* \cF \ra Q$ is zero.
Since this composite coincides with $g^*(\pi^* \cE^{\vee} \ra \pi^* \cF \ra \cQ)$,
we see that 
$g$ factors through $  \Gr (\cC,r) \hookrightarrow  \Gr (\cF,r) $ if and only if 
$g$ factors through the zero locus of $\pi^* \cE^{\vee} \ra \pi^* \cF \ra \cQ$.
\end{proof}

For the case $r = 1$, we get the projective bundle $\bP (\cF)$ and write $p_{\cF}$ for the projection of $\bP (\cF)$ to $M$. We can also view the composition of $p_{\cF}^{\ast} \sigma \colon p_{\cF}^{\ast}\cE^{\vee} \to p_{\cF}^{\ast}\cF$ with the canonical map $p_{\cF}^{\ast} \cF \to \cO_{\cF} (1)$ as a global section $s_{\sigma}$ of the bundle
\begin{equation} \label{Cayleytrick}
    \sHom (p_{\cF}^{\ast} \cE^{\vee}, \cO_{\cF} (1)) \cong \cE \boxtimes \cO_{\cF} (1).
\end{equation}
Here and subsequently, we use $\cV \boxtimes \cO_{\cF} (\ell)$ to denote $p_{\cF}^{\ast} \cV \otimes \cO_{\cF} (\ell)$ for any bundle $\cV$ over $M$ and $\ell \in \bZ$ by abuse of notation.

The following lemma is proved in \cite[Lemma 3.5 and 3.6]{LW2022}. 

\begin{lemma}[\cite{LW2022}]\label{detzer}
For $e \geqslant f$, we have:
\begin{enumerate}
    \item\label{detzer1} The projective bundle $\bP (\coker \sigma)$ coincides with the zero scheme $Z (s_{\sigma})$ in $\bP (\cF)$.
    
    \item\label{detzer2} The restriction of $p_{\cF}$ to $Z (s_{\sigma})$ maps onto $D_{f - 1}(\sigma)$ and it is an isomorphism if $D_{f - 2} (\sigma) = \varnothing$.
\end{enumerate}
\end{lemma}

Note that the expected codimension of $Z (s_{\sigma})$ in $\bP (\cF)$ is $e$, and given $x \in D_{f - 1} (\sigma)$ the fiber of $Z (s_{\sigma})$ over $x$ is $\bP(\coker \sigma (x))$. Also, Lemma \ref{detzer} \eqref{detzer1} is the case $r = 1$ in Lemma \ref{eq_as zero locus}.

\subsection{Calabi--Yau threefolds}

We recall a notion of the generality of morphisms of bundles used in \cite{SSW18}. Let $M$ be a smooth projective variety and let $\sigma \colon \cE^\vee \to \cF$ be a morphism of vector bundles on $M$ of rank $e$ and $f$ respectively.

\begin{definition}
For a given integer $r \geqslant 0$, the $\sigma$ is said to be $r$-general if $D_i (\sigma) \setminus D_{i - 1} (\sigma)$ is smooth of (expected) codimension $(e - i) (f - i)$ in $M$ for all $i = 0, \cdots, r$.
\end{definition}

In the remainder of this subsection we assume that $\cE$ and $\cF$ have the same rank $n + 1$ with the Calabi--Yau condition 
\begin{equation}\label{CYcond}
        \det (\cE) \otimes \det (\cF) \cong \cO (- K_M)
        \end{equation}
and $\dim M = 4$. To construct Calabi--Yau threefolds, we need the following two propositions, see \cite[Proposition 3.6, Theorem 4.4]{SSW18} and \cite[Remark 3.4, Proposition 4.4 (i) and B.1]{LW2022}.

\begin{proposition}[\cite{SSW18}]\label{n-general}
With the above notation, we fix an $n$-general morphism $\sigma$. Let us denote by $X_\cF$ the zero scheme $Z(s_{\sigma}) \subseteq \bP({\cF})$ and by $\pi_{\cF}$ 
the restriction of $p_{\cF}$ to $X_{\cF}$. Then $X_\cF$ is smooth and $\bP(\coker \sigma (x)) = \bP^1$ for $x \in D_{n - 1} (\sigma)$. 
\end{proposition}

Hence the dimension of the exceptional set $\pi_\cF$ is one, that is, the resolution is \emph{small}.

\begin{proposition}[\cite{LW2022}] \label{nodal3F}
With the above notation, we assume that the bundle $\cE \otimes \cF$ is ample and globally generated. Then there exists a Zariski open set $U$ in $\Hom (\cE^{\vee}, \cF)$ such that for each $\sigma \in U$, the morphism $\sigma$ is $n$-general and the $X_\cF = Z(s_{\sigma})$ is a Calabi--Yau threefold with Picard number $\rho (\bP (\cF))$. 
Moreover, the $D_n (\sigma) \subseteq M$ is a nodal hypersurface with the singular locus $D_{n - 1} (\sigma)$ and 
\begin{align*}
    \chi_{top} (X_\cF) &= \chi_{top} (\tY) + 2 |D_{n - 1} (\sigma)|, \\
      h^{2,1} (X_\cF) &= h^{2,1} (\tY) - |D_{n - 1} (\sigma)| + 1,
\end{align*}
where $\tY \in |- K_M|$ is a smooth member and $\chi_{top}(-)$ is the topological Euler number.
\end{proposition}

In Section \ref{V4_sec}, \ref{V5_sec} and \ref{G24_sec}, we will apply Proposition \ref{nodal3F} with $\cE = \cO^{\oplus (n + 1)}$ and Fano bundles $\cF$ listed in Theorem \ref{clFano}. We remark that, by the classification, such a Fano bundle $\cF$ is globally generated (cf.\ the proof of \cite[Proposition 2.7]{LW2022}).

The $\pi_{\cF} \colon X_{\cF} \to D_n (\sigma)$ is called the \emph{determinantal contraction} of $X_{\cF}$. The singular locus $D_{n - 1} (\sigma)$ consists of ordinary double points (ODPs for short). 

The following is from \cite[Remark 3.3, Proposition 4.5]{SSW18}. It computes some invariants of $X_{\cF}$ in terms of Chern classes of the virtual bundle $\cE - \cF^{\vee}$ and its dual. For the formulas for Chern classes of such bundles, we refer the reader to \cite[Example 3.2.7]{fulton98} or \cite[Appendix A]{LW2022}.

\begin{proposition}[\cite{SSW18}] \label{LHinterODP}
With the above notation, we fix an $n$-general morphism $\sigma$. Let $H_M$ be any Cartier divisor on $M$, $H_\cF = \pi_{\cF}^*H_M$ and $L_\cF = \cO_{\cF} (1)|_{X_\cF}$. Then we have
\begin{align*}
    \int_{X_{\cF}} H_{\cF}^k \cdot L_{\cF}^{3-k} &= \int_M H_M^k\cdot c_{4 - k}(\cE - \cF^{\vee}) \text{ for } 0 \leqslant k \leqslant 3,\\
    \int_{X_{\cF}} c_2(T_{X_{\cF}})\cdot H_{\cF} &= \int_M c_2(T_M) \cdot c_1(\cE - \cF^{\vee})\cdot H_M, \\
    \int_{X_{\cF}} c_2(T_{X_{\cF}})\cdot L_{\cF} &= \int_M c_2(T_M)\cdot c_2(\cE - \cF^{\vee}) - |\Sing (D_n (\sigma))|,
\end{align*}
and the number of ODPs of $D_{n} (\sigma)$ is
\[
    \int_M c_2 (\cF -\cE^{\vee})^2 - c_1 (\cF -\cE^{\vee}) \cdot c_3 (\cF -\cE^{\vee}).
\]
\end{proposition}

There is the other determinantal contraction $\pi_{\cE} \colon X_{\cE} \to D_{n} (\sigma^{\vee})$ via the dual morphism $\sigma^{\vee} \colon \cF^{\vee} \to \cE$. By the observation that $D_n (\sigma) = D_n (\sigma^{\vee})$, we get a birational map $\chi \coloneqq \pi_{\cE}^{- 1} \circ \pi_{\cF}$. The following lemma gives some information of the proper transform $\chi_\ast L_\cF$, see \cite[Lemma 4.8]{LW2022} for a proof.

\begin{lemma}\label{ab}
Under the assumptions as in Proposition \ref{LHinterODP}, we assume that $\chi_* L_{\cF} = \alpha D_1 + \beta D_2$ in $\Pic(X_{\cE})_\bQ$ and $D_1$, $D_2$ are base point free. If $\chi$ is not an isomorphism, then $\alpha \beta < 0$,
\[
    L_{\cF} \cdot H_{\cF}^2 = \chi_*L_{\cF} \cdot (\chi_* H_{\cF})^2 \text{ and } L_{\cF}^2 \cdot H_{\cF} = (\chi_*L_{\cF})^2\cdot \chi_* H_\cF.
\]
\end{lemma}

\subsection{Flops via dual morphisms}\label{Nefcone_subsec}

Let $(M, \cF)$ be a pair as in Theorem \ref{clFano}, where $\rank \cF = n + 1$, and $n =1$ or $2$. Let  $\cE = \cO_M^{\oplus n + 1}$. We shall apply Proposition \ref{nodal3F} to construct our Calabi--Yau threefolds, which will be used in the following sections. It gives rise to smooth Calabi--Yau threefolds $X_\cF$ and $X_\cE$ with Picard number two for a general morphism $\sigma \colon \cE^{\vee} \to \cF$. Since $\cE$ is trivial, the Calabi--Yau $X_{\cF} \subseteq \bP (\cF)$ is a complete intersection of $n + 1$ general hypersurfaces in $|\cO_\cF (1)|$. 

The determinantal contractions and the restriction of the projection $\bP (\cE) = M \times \bP^n \to \bP^n$ to $X_\cE$ give the following diagram
\begin{equation} \label{flopdiag}
    \begin{tikzcd}
    X_\cF \ar[dr,"\pi_{\cF}"] \ar[rr,dashed,"\chi"] & & X_{\cE} \ar[dl,"\pi_{\cE}", swap] \ar[d]  \\
    & D_n(\sigma) & \bP^n.  \\
    \end{tikzcd}
\end{equation}
We assume that there exists another small contraction $X_\cF \to Y_\cF$.

We let $H_M$ be the fundamental divisor of the Fano fourfold $M$, and write $L_\cF$ and $H_\cF$ for the restrictions of $\cO_\cF (1)$ and $p_{\cF}^\ast H_M$ to the Calabi--Yau $X_\cF$, and similarly for $L_\cE$ and $H_\cE$.

Recall that $n = 1$ if $M  = V_4$ or $V_5$, and $n =2$ if $M = \Gr (2, 4)$. According to $\dim X_{\cE} =3 > n$, it follows that the fibration $X_\cE \to \bP^n$ gives one edge of the nef cone $\Nef(X_\cE)$ of $X_\cE$. Since $X_\cE$ has Picard number $2$, the diagram \eqref{flopdiag} implies that
$\Nef(X_\cE)$ is spanned by $L_{\cE}$ and $H_{\cE}$. To describe the cone in $N^1 (X_\cF)$, we need to compute some intersection numbers. For abbreviation, we write $L$ and $H$ instead of $L_\cF$ and $H_\cF$ and derive Table \ref{intAll} by using Proposition \ref{LHinterODP}.

\begin{table}
  \centering
  \begin{tabular}{cccccccc}
    \toprule
    $(M, \cF)$ & $L^3$ & $L^2.H$ & $L.H^2$ & $H^3$ & $L.c_2(T_{X_{\cF}})$ & $H.c_2(T_{X_{\cF}})$ & $\#$ of ODPs \\
    \midrule
    $(V_4, p^{\ast} \cS^{\vee} (1))$ & 80 & 48 & 26 & 12 & 104 & 60 & 26 \\
    $(V_5, \cS_{V_5}^{\vee} (1))$ & 110 & 63 & 33 & 15 & 116 & 66 & 29\\
    $(\Gr(2, 4), \cS (2) \oplus \cO (1))$ & 85 & 45 & 21 & 8 & 106 & 56 & 41\\
    \bottomrule
  \end{tabular}
  \caption{The intersection numbers of $X_\cF$.}
  \label{intAll}
\end{table}

\begin{proposition}\label{movXE}
Under the above assumptions, the matrix representation of the map $\chi_{\ast} \colon N^1(X_\cF)\ra N^1(X_\cE)$ with respect to $\{L,H\}$ and $\{L_\cE,H_\cE\}$ is given by
\[
    [\chi_*] =
    \begin{bmatrix}
      -1 & 0 \\
      r_M & 1
    \end{bmatrix}
    =
    [(\chi^{-1})_*].
\]
where $r_M$ is the index of the smooth Fano fourfold $M$. Moreover, in $N^1 (X_{\cF})$ we have 
\[
    \Nef (  X_{\cE}) = \bR_{ \geqslant 0} [r_M H - L]+ \bR_{ \geqslant 0} [H].
\]
\end{proposition}

\begin{proof}
We first claim that $\chi$ is not an isomorphism. If the claim is not true, then $X_\cF \cong X_\cE$ admits a fibration over $\bP^n$, where $n = 1$ or $2$. This is contrary to the assumption that $X_\cF$ admits another small contraction $X_\cF \to Y_\cF$.  

Write $\chi_*L = \alpha L_\cE + \beta H_\cE$. Since $\chi$ is an isomorphism in codimension one, the proper transform $\chi_{\ast} H$ is $H_{\cE}$. By Proposition \ref{LHinterODP}, we get that $L_\cE \cdot H_{\cE}^2 = 10, 12, 11$ and $L_{\cE}^2 \cdot H_{\cE} = 0, 0 , 5$ for $M  = V_4, V_5, \Gr (2, 4)$ respectively. From Table \ref{intAll} and Lemma \ref{ab}, we get a system of equations for $\alpha$ and $\beta$. Solving such system yields $\alpha = \pm 1$. According to $\alpha \beta < 0$, it follows that $(\alpha, \beta) = (-1, r_M)$. Note that $r_{V_4} = r_{V_5} = 3$ and $r_{\Gr(2,4)} = 4$. 

Since the projection $X_\cE \to \bP^n$ and the determinantal contraction $\pi_\cE$ are induced by $L_\cE$ and $H_\cE$ respectively, we find that $\bR_{ \geqslant 0} [r_M H - L]$ and $\bR_{ \geqslant 0} [H]$ are the two boundary rays of $\Nef (X_\cE)$ in $N^1 (X_\cF)$, which completes the proof.
\end{proof}

Finally, we deduce the Hodge numbers of $X_\cF$ in Table \ref{topAll} by applying Proposition \ref{nodal3F}, Table \ref{intAll} and \cite[Table 6]{LW2022}.

\begin{table}
  \centering
  \begin{tabular}{cccccccc}
    \toprule
    $(M, \cF)$ & $\chi_{top} (X_\cF)$ & $h^{2,1} (X_\cF)$ \\
    \midrule
    $(V_4, p^{\ast} \cS^{\vee} (1))$ & -92 & 48 \\
    $(V_5, \cS_{V_5}^{\vee} (1))$ & -92 & 48 \\
    $(\Gr(2, 4), \cS (2) \oplus \cO (1))$ & -94 & 49 \\
    \bottomrule
  \end{tabular}
  \caption{The topological numbers of $X_\cF$.}
  \label{topAll}
\end{table}

\section{The Del Pezzo fourfold of degree 4}\label{V4_sec}

We consider the case $(M, \cF, \cE) = (V_4, p^{\ast} \cS^{\vee} (1), \cO^{\oplus 2})$ and will use the same notation as in Section \ref{Nefcone_subsec}. For a general morphism $\sigma \colon \cE^\vee \to \cF$,  there are two smooth Calabi-Yau threefolds $X_{\cF}$, $X_{\cE}$ with the determinantal  contractions $\pi_{\cF}$, $\pi_{\cE}$ and the diagram \eqref{flopdiag}.

Recall that  $p \colon V_4 \to \Gr(2,4)$ is a double cover branching along a smooth divisor $B \in |\cO_{\Gr(2,4)}(2)|$.
Furthermore, the $\cS $ is the universal subbundle on $\Gr(2,4)$.
Hence we have the following diagram:
\begin{equation*} \label{diagram_V4}
    \begin{tikzcd}
     & \bP_{V_4}(\cF) \ar[d,"\tilde{p}"] \ar[dr,"p_{\cF}"] \ar[ddl,bend right,"q"'] & \\
     & \bP_{\Gr(2,4)}(\cS^{\vee}) \ar[dr,"p_2"] \ar[dl, "p_1"'] & V_4 \ar[d,"p"]  \\
    \bP^3 &  &  \Gr(2,4).  
    \end{tikzcd}
\end{equation*}
Here $\tilde{p} $ is the double cover branching along $ p_2^{-1}(B)$
and $ p_1$ is the morphism defined by $|\cO_{\cS^{\vee}}(1)|$.
We note that $\bP_{\Gr(2,4)}(\cS^{\vee}) $ is the flag variety $F(1,2;4)$ and $p_1$ is nothing but the natural morphism 
$ F(1,2;4) \to \Gr(1,4)= \bP^3$.
The morphism $q$ is just the composite $p_1 \circ \tilde{p} $. 

Since $\cF =p^* \cS^{\vee} (1)$, we have $q^{\ast} \cO_{\bP^3} (1) \cong \cO_\cF (1) \boxtimes \cO_{V_4} (- 1)$.
Consider the restriction of $q$ to $X_{\cF}$ and let $X_{\cF} \to Y_{\cF} \to \bP^3$ be the Stein factorization. 
In this section, we prove the following theorem:

\begin{theorem}\label{mainV4}
Let $(M, \cF, \cE) = (V_4, p^{\ast} \cS^{\vee} (1), \cO^{\oplus 2})$. Then for a general morphism $\sigma \colon \cE^{\vee} \to \cF$, the scheme $X_{\cF}$ is a smooth Calabi--Yau threefold of Picard number two with
\[
    \Nef(X_\cF)=\bR_{\geqslant 0}[L-H]+\bR_{\geqslant 0}[H],
\]
such that
\begin{enumerate}[(i)]
    \item\label{mainV41} the determinantal contraction $\pi_{\cF}$ is induced by $H$;

    \item\label{mainV42} the morphism $X_{\cF} \to \bP^3 $ is generically $2:1$,
    which induces an involution $\iota \colon X_{\cF} \dashrightarrow X_{\cF}$ over $\bP^3$;
    
    \item\label{mainV42.5} the stein factorization $X_{\cF} \to Y_{\cF}$ of $X_{\cF} \to \bP^3 $ is a small contraction induced by $L - H$;
    

    \item\label{mainV43} the $X_{\cE}$ admits a K3 fibration induced by $3 H - L$.
\end{enumerate}
Moreover, the movable cone $\cMov(X_{\cF})$ is the convex cone generated by the divisors $15 L - 17 H$ and $3 H - L$ which is covered by the nef cones of $X_{\cF}$ and $X_{\cE}$, and there are no more minimal models of $X_{\cF}$, which we summarize in the following diagram:
\begin{equation*}
   \begin{tikzcd}[row sep=small]
      X_{\cE} \ar[d] \ar[rd, "\pi_{\cE}"] \ar[rr, dashrightarrow, "\chi^{-1}"] && X_{\cF} \ar[ld, "\pi_{\cF}", swap] \ar[rdd] \ar[d] \ar[rr, dashrightarrow, "\iota"] & & X_{\cF} \ar[rr, dashrightarrow, "\chi"] \ar[ldd] \ar[d] \ar[rd, "\pi_{\cF}"] && X_{\cE} \ar[ld, "\pi_{\cE}", swap] \ar[d] \\
     \bP^1& D_1 (\sigma) & Y_{\cF} \ar[rd] \ar[rr, crossing over] & &Y_{\cF} \ar[ld] & D_1 (\sigma) & \bP^1 . \\
     &&& \bP^3 &&&
    \end{tikzcd}
\end{equation*}
\end{theorem}

Figure \ref{movV4sliPic} is the slice of the movable cone of $X_\cF$. It is a subdivision of a closed interval, which comes from the chamber structure of the cone. We depict $X_{\cF}$ and $X_{\cE}$ inside their nef cones. 

\begin{figure}[htbp]
\centering
\begin{tikzpicture}
\draw(0,0)--(10,0);

\foreach \x in {0,2.5,5,7.5,10}
  \draw[very thick, teal] (\x,5pt)--(\x,-5pt);

\foreach \y/\ytext in {-0.25/$15L - 17H$,2.5/$8 L - 9 H$,5/$L - H$,7.5/$H$,10.25/$3 H - L$}
  \draw (\y,0) node[above=1ex] {\ytext};

\foreach \z/\ztext in {1.25/$X_{\cE}$,3.75/$X_{\cF}$,6.25/$X_{\cF}$,8.75/$X_{\cE}$}
  \draw (\z,0) node[below] {\ztext};
\end{tikzpicture}
\caption{The slice of the movable cone $\cMov(X_\cF)$ for $M= V_4$.} \label{movV4sliPic}
\end{figure}






\vspace{2mm}
In the rest of this section, we prove Theorem \ref{mainV4}.

\begin{lemma}\label{V4small}
If $\sigma \colon \cE^{\vee} \to \cF$ is chosen in general, then the natural projection $X_{\cF} \to \bP^3$ is generically $2:1$ and the stein factorization $X_{\cF} \to Y_{\cF}$ of $X_{\cF}  \to \bP^3$ is a small contraction induced by $L - H$.
\end{lemma}

\begin{proof}
Since $\cO_{\cF} (1) \boxtimes \cO_{V_4} (- 1) \cong q^{\ast} \cO_{\bP^3} (1)$, it holds that $L - H \cong q|_{X_{\cF}}^\ast \cO_{\bP^3}(1)$. Hence the morphism $X_{\cF} \to Y_{\cF}$ is induced by $L - H$.

Let $\cV = \Omega_{\bP^3} (2) \oplus \cO_{\bP^3}$. By \cite[Theorem 1.3 (4)]{NO07} (also see p.236, Proof of Theorem 1.3, Case 2) and \cite[Remark 4.2]{Kan19}, there is a divisor $Q \in |\cO_{\cV} (2)|$ such that 
\begin{align*}
    \bP_{V_4} (\cF) \cong Q \subseteq \bP_{\bP^3} (\cV).
\end{align*}
According to $c_1 (\Omega_{\bP^3} (2) \oplus \cO_{\bP^3}) = c_1 (\cO_{\bP^3} (2))$, it follows that 
\[
    \cO (K_{\bP (\cV)}) \cong \cO_\cV (- 4) \boxtimes \cO_{\bP^3} (- 2)
    \quad \text{and} \quad 
    \cO (K_{Q}) \cong \cO_\cV (- 2) \boxtimes \cO_{\bP^3} (- 2)|_Q.
\]
Since $\cO_{\cF}(-2) =\cO(K_{ \bP_{V_4} (\cF)}) = \cO(K_Q)$ under the identification $\bP_{V_4} (\cF) \cong Q$,
we have $\cO_{\cF}(1) \cong \cO_\cV (1) \boxtimes \cO_{\bP^3} (1)|_Q$ as the Picard group of a Fano manifold is torsion free. 
Note that $\cV$ is globally generated because we have 
\[
    0 \to \Omega_{\bP^3}^2 (2) \to \cO_{\bP^3}^{\oplus 6} \to \Omega_{\bP^3} (2) \to  0   
\]
by dualizing the Euler sequence and taking the $2$nd exterior power. Hence $|\cO_\cV (1) \boxtimes \cO_{\bP^3} (1))|$ is ample and base point free (see \cite[Lemma 2.6]{LW2022}).
Furthermore, the natural map $H^0 (\cO_\cV (1) \boxtimes \cO_{\bP^3} (1)) \to H^0 (\cO_\cV (1) \boxtimes \cO_{\bP^3} (1)|_Q)$ is surjective by 
\[
    H^1 (\cO_\cV (- 1) \boxtimes \cO_{\bP^3} (1)) = 0,
\] 
which follows from Kodaira vanishing theorem 
since 
\[
    \cO_\cV (- 1) \boxtimes \cO_{\bP^3} (1) = \cO (K_{\bP (\cV)}) \otimes (\cO_\cV (3) \boxtimes \cO_{\bP^3} (3)). 
\]
Hence the Calabi--Yau threefold $X_{\cF} \subseteq \bP_{V_4} (\cF) \cong Q$ is isomorphic to a complete intersection $D_1 \cap D_2 \cap Q$ where $D_1$ and $D_2$ are general members of $|\cO_\cV (1) \boxtimes \cO_{\bP^3} (1)|$. 

Let $V \subseteq \bP_{\bP^3} (\cV)$ be the complete intersection fourfold defined by $D_1$ and $D_2$, i.e., 
it is the zero locus of $ U \otimes \cO_{\bP_{\bP^3} (\cV)} \to \cO_\cV (1) \boxtimes \cO_{\bP^3} (1)$
induced by a general two-dimensional subspace $U$ of $H^0(\cO_\cV (1) \boxtimes \cO_{\bP^3} (1))$. 
Let $\tilde{q} \colon \bP_{\bP^3}(\cV) \to \bP^3 $ be the natural projection. Then $U$ also defines a morphism $\tau \colon U \otimes \cO_{\bP^3} \to \cV \otimes \cO_{\bP^3} (1)$ by the natural isomorphism $\tilde{q}_{\ast} \cO_{\cV} (1) \cong \cV$.

Consider the restriction $q_V \colon V \to \bP^3$ of $\tilde{q}$ on $V$. For each $x \in D_k (\tau) \setminus D_{k - 1} (\tau)$, the fiber $q_V^{- 1} (x)$ is $\bP (\coker \tau (x)) \cong \bP^{4 - k - 1}=\bP^{3-k}$ for $0 \leqslant k \leqslant 2$ by Lemma \ref{detzer}. Note that $D_1 (\tau) $ and $D_0 (\tau)$ have the expected codimension $(2 - 1) (4 - 1) = 3$ and $2 \cdot 4 = 8$ respectively. If $U$ is chosen in general, then $D_0 (\tau) = \varnothing$ and $D_1(\tau)$ consists of finitely many (smooth) points whose number can be determined by Giambelli--Thom--Porteous formula \cite[Theorem 14.4]{fulton98},
\[
    [D_1(\tau)] = c_3 (\cV \otimes \cO_{\bP^3} (1)) \cap [\bP^3] \in A_3(\bP^3).
\]
According to $c_3 (\cV \otimes \cO_{\bP^3} (1)) = 14$, the projection $q_V \colon V \to \bP^3$ has exactly 14 projective planes $\bP^2$ as fibers, and the other fibers of $q_V$ are $\bP^1$. Then $q_V^{-1} (x) \cap Q$ is a quadric in $\bP^2$ or a $\bP^1$ or a (possibly) double point in $\bP^1$. 
In particular, $X_{\cF} = V \cap Q \to \bP^3$ is generically $2:1$.
We note that this also follows from $(L - H)^3 = 2$ computed by Table \ref{intAll}.

\begin{claim}
\label{claim_small}
The exceptional set of $X_{\cF} \to Y_{\cF}$ has dimension at most one.
\end{claim}
\begin{proof}[Proof of Claim \ref{claim_small}]
We denote $V$, $\tau$ by $V_U$, $\tau_U$ when we need to clarify the choice of $U$.  
By \cite[Theorem 4.7]{Kan19}, the morphism $q \colon Q \to \bP^3$ is a quadric bundle,
that is, the fiber
$q^{-1}(x) = Q \cap \tilde{q}^{-1}(x)$ is a quadric in $\tilde{q}^{-1}(x) \cong \bP^3 $ for any $x \in \bP^3$.
Let $\Gr$ be the Grassmannian of two-dimensional subspaces in $H^0(\cO_{\cV}(1) \boxtimes \cO_{\bP^3}(1))$ and $\Sigma \subseteq \bP^3$ the discriminant locus of $q \colon Q \to \bP^3$. 
We set
\[
    \Gr^{\circ} \coloneqq \{[U] \in \Gr \mid  \text{$ \dim D_1(\tau_U) =0$, $D_1(\tau_U ) \cap \Sigma = D_0(\tau_U) =\varnothing$} \}.
\]
This is a non-empty open subset of $\Gr$. Indeed, we already saw that $ \dim D_1(\tau_U) =0$ and $D_0(\tau_U) =\varnothing$ are satisfied if $U$ is general.
Since $\dim \Sigma=2$ and the expected codimension of $D_1(\tau|_{\Sigma})$ in $\Sigma$ is three,
we have $D_1(\tau|_{\Sigma}) = D_1(\tau ) \cap \Sigma  =\varnothing$ if $U$ is general.

For each $x \in \bP^3$, the natural map 
\begin{align}\label{eq_restriction}
 H^0(\cO_{\cV}(1)\boxtimes \cO_{\bP^3}(1)) \rightarrow
H^0(\cO_{\cV}(1)\boxtimes \cO_{\bP^3}(1)|_{\tilde{q}^{-1}(x) } )  \cong H^0(\cO_{\bP^3}(1))
\end{align}
is surjective.
This is because this map can be identified with 
$  H^0(\cV\otimes \cO_{\bP^3}(1)) \rightarrow  \cV\otimes \cO_{\bP^3}(1) \otimes k(x)$,
which is surjective since $\cV$ is globally generated.
Furthermore,
the restriction map 
\begin{align*}\label{eq_restriction2}
H^0(\cO_{\cV}(1)\boxtimes \cO_{\bP^3}(1)|_{\tilde{q}^{-1}(x) } ) 
\to H^0(\cO_{\cV}(1)\boxtimes \cO_{\bP^3}(1)|_{q^{-1}(x) } ) 
\end{align*}
is an isomorphism for any $x \in \bP^3$ since $q^{-1}(x) $ is a quadric in $\tilde{q}^{-1}(x) \cong \bP^3 $.

To prove Claim \ref{claim_small},
it suffices to show that $X_{\cF} = Q \cap V_U \to \bP^3$ contracts at most finitely many curves for general $U$.
Observe that the fiber of $Q \cap V_U \to \bP^3$ over $x \in \bP^3$ is 
\[
q^{-1}(x) \cap q_{V_U}^{-1}(x) \subseteq \tilde{q}^{-1}(x) \cong \bP^3.
\]
If $[U] \in \Gr^{\circ}$,
we have $\dim q_{V_U}^{-1}(x) \cap  q^{-1}(x)  \leqslant 1$ for any $x \in \bP^3$.
In fact, this inequality is trivial if $x \not \in  D_1(\tau_U)$ since  $q_{V_U}^{-1}(x) \cong \bP^1 $ in this case.
If $x \in D_1(\tau_U) $, 
the $q_{V_U}^{-1}(x) \cong \bP^2$ is a hyperplane of $\tilde{q}^{-1}(x) $.
By $D_1(\tau_U) \cap \Sigma =\varnothing$ and $x \in D_1(\tau_U)$,
the $q^{-1}(x)$ is a smooth quadric and hence the hyperplane section $q^{-1}(x) \cap q_{V_U}^{-1}(x) $ is one-dimensional for $x \in D_1(\tau_U)$.

Thus the rest is to show that $q^{-1}(x) \cap q_{V_U}^{-1}(x)$ is one-dimensional for at most finitely many $x \in \bP^3$
if $U$ is general.
To see this, we consider 
\[
Z \coloneqq \{ ([U], x) \in \Gr^{\circ} \times \bP^3 \mid \dim q^{-1}(x) \cap q_{V_U}^{-1}(x) = 1\}.
\]
For $[U] \in \Gr^{\circ}$,
the morphism $X_{\cF} = Q \cap V_U \to \bP^3$ contracts a curve over $x \in \bP^3$ 
if and only if $ ([U],x) $ is contained in $Z$.
Thus if $\dim Z \leqslant \dim \Gr^{\circ}$, the fiber of $Z \to \Gr^{\circ}$ over a general $[U]$ is at most zero-dimensional
and we obtain the finiteness of $x \in \bP^3$ with $\dim q^{-1}(x) \cap q_{V_U}^{-1}(x) =1$.

We are left with the task to check $\dim Z \leqslant \dim \Gr^{\circ}$.
For $x \in \bP^3$, let 
\[
    K_x \subseteq  H^0(\cO_{\cV}(1)\boxtimes \cO_{\bP^3}(1))
\] 
be the kernel of \eqref{eq_restriction},
whose codimension is $\dim H^0(\cO_{\cV}(1)\boxtimes \cO_{\bP^3}(1)|_{\tilde{q}^{-1}(x) } ) =4$.
Hence the codimension of the Schubert variety
\begin{align*}
  B_x \coloneqq \{[U] \in \Gr \mid  U \cap K_x \neq \{0\} \} = \{[U] \in \Gr \mid  x \in D_1(\tau_U)\} 
\end{align*}
in $\Gr$ is three.
For $[U] \in \Gr^{\circ} \setminus B_x$, the $q_{V_U}^{-1}(x)$ is a line since $x \not \in D_1(\tau_U)$.
Hence we can define the morphism 
\begin{equation}\label{eq_Gr_restriction_x}
 \Gr^{\circ} \setminus B_x \to \Gr(2,H^0(\cO_{\cV}(1)\boxtimes \cO_{\bP^3}(1)|_{\tilde{q}^{-1}(x) } )  ) =\{ \text{lines in $\tilde{q}^{-1}(x)$}\}
\end{equation}
by sending $[U] \in \Gr^{\circ} \setminus B_x$ to the line $q_{V_U}^{-1}(x)$ in $\tilde{q}^{-1}(x) \cong \bP^3$.
If $x \not \in \Sigma$, the family of lines in the smooth quadric $q^{-1}(x) \cong \bP^1 \times \bP^1$ is one-dimensional and hence
the codimension of $\{ \text{lines in $q^{-1}(x)$}\} $ in $\{ \text{lines in $\tilde{q}^{-1}(x)$}\} \cong \Gr(2,4)$ is three.
Similarly,
the codimension of $\{ \text{lines in $q^{-1}(x)$}\} $ in $\{ \text{lines in $\tilde{q}^{-1}(x)$}\}$ is two or three
if $x \in \Sigma$.
Thus $\{ [U] \in \Gr^{\circ} \setminus B_x \, | \, q_{V_U}^{-1}(x) \subseteq q^{-1}(x) \}$, if it is non-empty, is of codimension three in $\Gr^{\circ} \setminus B_x$ when $x \not \in \Sigma$,
and of codimension at least two when $x \in \Sigma$.
Note that the morpshism \eqref{eq_Gr_restriction_x} is dominant and all the fibers are isomorphic to open subsets of $\Gr(2, h^0(\cO_{\cV}(1) \boxtimes \cO_{\bP^3}(1))-2)$.
Hence the codimension does not change after the pullback by \eqref{eq_Gr_restriction_x}.

By definition,
the fiber of $Z \to \bP^3$ over $ x \in \bP^3$
is contained in 
\[
\{ [U] \in \Gr^{\circ}\setminus B_x \, | \, q_{V_U}^{-1}(x) \subseteq q^{-1}(x) \} \cup (B_x \cap \Gr^{\circ}),
\]
whose codimension in $\Gr^{\circ}$ is three (resp.\ at least two) if $x \not \in \Sigma $ (resp.\ $x \in \Sigma$),
if it is non-empty.
Since $x$ moves on $\bP^3$ and $\dim \Sigma =2$,
the dimension of $Z$ is at most  that of $\Gr^{\circ}$, and Claim \ref{claim_small} follows.
\end{proof}

We note that $X_{\cF} \to Y_{\cF}$ is not an isomorphism. Indeed, we have already seen that $X_\cF$ is smooth, irreducible and has Picard number two. Notice that $Y_\cF \to \bP^3$ is a finite morphism of degree two. If $X_{\cF} \to Y_{\cF}$ is an isomorphism, then the Picard number of $X_{\cF} \cong Y_\cF$ coincides with that of $\bP^3$ (see, e.g., \cite[Example 7.1.20]{LazPAGII}), a contradiction.
Hence $X_{\cF} \to Y_{\cF}$ is a small contraction by Claim \ref{claim_small}. 
\end{proof}

\begin{remark}
Recall that $B \in |\cO_{\Gr(2,4)}(2)|$ is the branched divisor of the double cover $p \colon V_4 \to \Gr(2,4)$.
If we assume that $B$ is general, so is the $Q \in |\cO_\cV (2)|$.
In this case, we can check Claim \ref{claim_small} more easily, and show that $Y$ has only $102$ ODPs.
\end{remark}


There is an involution  on $X_\cF$ over  $\bP^3$ induced from the generically $2:1$ morphism $X_{\cF} \to \bP^3$. We denote the involution by $\iota \colon X_\cF \dra X_\cF$.

\begin{lemma} \label{V4invol}
For the involution $\iota \colon X_\cF\dra X_\cF$, the matrix representation
of $\iota^{\ast} : N^1(X_{\cF}) \to N^1(X_{\cF}) $ with respect to $\{L,H\}$ are given by
\[
\begin{bmatrix}
\iota^\ast
\end{bmatrix}
=
\begin{bmatrix}
9 & 8 \\
-10 & -9
\end{bmatrix}
=
\begin{bmatrix}
(\iota^{-1})^\ast
\end{bmatrix}
\]
\end{lemma}

\begin{proof} 
Note that $\iota^*(L-H)=L-H$ 
since $L-H$ is the pullpack of $\cO_{\bP^3}(1)$ by $X_{\cF} \ra \bP^3$. 
We let $\bar{q} = q|_{X_\cF} \colon X_\cF \to \bP^3$ and $\bar{q}_* H = a H_{\bP^3}$, where $a \in \bZ$ and $H_{\bP^3}$ is the hyperplane class of $\bP^3$. Then
\[
    a = (\bar{q}_\ast H \cdot H_{\bP^3}^2)_{\bP^3} = (H \cdot (L - H)^2)_{X_\cF} = 8
\]
by the projection formula and Table \ref{intAll}. Hence 
\[
    H + \iota^\ast H = \bar{q}^\ast \bar{q}_\ast H = \bar{q}^\ast (8 H_{\bP^3}) = 8 (L - H).
\]
Then $\iota^\ast H = 8 L - 9 H$ and the rest is clear.
\end{proof}


We are now in a position to prove Theorem \ref{mainV4}.


\begin{proof}[Proof of Theorem \ref{mainV4}]
The \eqref{mainV41} follows from the construction of $\pi_\cF$, and \eqref{mainV42}, \eqref{mainV42.5} from Lemma \ref{V4small}. 
By Proposition \ref{movXE},
the morphism $X_{\cE} \to \bP^1$ is induced by $3L-H$.
Furthermore, $X_{\cE} \to \bP^1$ is a K3 fibration. Indeed, by Proposition \ref{LHinterODP}, we have
\begin{equation}\label{K324}
    \int_{X_{\cE}} c_2 (T_{X_{\cE}}) \cdot L_{\cE} = \int_{V_4} c_2 (T_{V_4}) \cdot c_2 (\cF) - c_2 (\cF)^2 = 24
\end{equation}
and then the intersection number of $c_2 (T_{X_{\cE}})$ with the general fiber $F \in |L_{\cE}|$ of $X_{\cE} \to \bP^1$ is $24$. Hence $F$ is a K3 surface by \cite[Lemma 3.3]{Oguiso93}, which proves \eqref{mainV43}.

Finally, we compute the boundaries of $\cMov(X_\cF)$. By Lemma \ref{V4invol} and Proposition \ref{movXE}, we find that $\iota_\ast H = 8 L - 9 H$ and $\iota_\ast \chi_\ast L_\cE = 15 L - 17H$, which completes the proof.
\end{proof}

\section{The Del Pezzo fourfold of degree 5}\label{V5_sec}

We consider the case $(M, \cF, \cE) = (V_5, \cS_{V_5}^{\vee} (1), \cO^{\oplus 2})$. 
Recall that $V_5$ is a general linear section of the Grassmannian $\Gr(2, 5) =\Gr(2,W) \subseteq \bP(\wedge^3 W)$ of codimension two, 
and $\cS_{V_5}$ is the restriction of the universal subbundle on $\Gr (2, W)$ to $V_5$,
where $ W$ is a five-dimensional vector space over $\bC$.  
In Section \ref{Nefcone_subsec}, we have already seen that for a general morphism $\sigma \colon \cE^\vee \to \cF$, there are two smooth Calabi-Yau threefolds $X_{\cF}$, $X_{\cE}$ and the diagram \eqref{flopdiag}.
In this section,
we prove the following theorem.


\begin{theorem}\label{mainV5}
Let $(M, \cF, \cE) = (V_5, \cS_{V_5}^{\vee} (1), \cO^{\oplus 2})$. Then for a general morphism $\sigma \colon \cE^{\vee} \to \cF$, the scheme $X_{\cF}$ is a smooth Calabi--Yau threefold of Picard number two with
\[
    \Nef(X_\cF)=\bR_{\geqslant 0}[L-H]+\bR_{\geqslant 0}[H],
\]
such that
\begin{enumerate}[(i)]
    \item\label{mainV51} the determinantal contraction $\pi_{\cF}$ is induced by $H$;

    \item\label{mainV52} the $L -  H$ defines a small contraction $X_{\cF} \to Y_{\cF}$ and $Y_{\cF} \subseteq \bP (W^{\vee})$ is a quintic hypersurface with $54$ ODPs;

    \item\label{mainV53} the flop $X_\cF^+$ of $X_{\cF} \to Y_{\cF}$ admits a small contraction $\psi \colon X_\cF^+ \to Z_{\cF}$ induced by $9 L - 11 H$ such that the exceptional locus $\Exc (\psi) $ is a smooth rational curve $ \Sigma^+$ 
    contracted to an ODP;
    
    \item\label{mainV54} the flop $X_\cF^{++}$ of $X_{\cF}^+ \to Z_{\cF}$ admits a K3 fibration induced by $4 L - 5 H$;
    
    \item\label{mainV55} $X_{\cE}$ admits a K3 fibration induced by $3 H - L$.
\end{enumerate}
Moreover, the movable cone $\cMov(X_{\cF})$ is the convex cone generated by the divisors $4 L - 5 H$ and $3 H - L$ which is covered by the nef cones of $X_{\cE}$, $X_{\cF}$, $X_{\cF}^+$ and $X_{\cF}^{++}$, and there are no more minimal models of $X_{\cF}$, which we summarize in the following diagram:
\begin{equation*}
    \begin{tikzcd}
       X_{\cF}^{++} \ar[d] \ar[rd] \ar[rr, dashrightarrow] && X_{\cF}^+ \ar[ld, "\psi"] \ar[rd, "\pi_\cH"] \ar[rr, dashrightarrow, "\theta"] && X_{\cF} \ar[ld, "\pi_{\wedge^3 \cT}"] \ar[rd, "\pi_{\cF}"] \ar[rr, dashrightarrow, "\chi"] && X_{\cE} \ar[ld, "\pi_{\cE}"] \ar[d]  \\
      \bP^1 & Z_{\cF} && Y_{\cF} && D_1 (\sigma) & \bP^1.
    \end{tikzcd}
\end{equation*}
\end{theorem}

The slice of the chamber structure of $\cMov(X_\cF)$ is given in Figure \ref{movV5sliPic}. We depict $X_{\cE}$, $X_{\cF}$, $X_{\cF}^+$ and $X_{\cF}^{++}$ inside their nef cones.

\begin{figure}[htbp]
\centering
\begin{tikzpicture}
\draw(0,0)--(10,0);

\foreach \x in {0,2.5,5,7.5,10}
  \draw[very thick, teal] (\x,5pt)--(\x,-5pt);

\foreach \y/\ytext in {-0.25/$4L - 5H$,2.5/$9 L - 11 H$,5/$L - H$,7.5/$H$,10.25/$3 H - L$}
  \draw (\y,0) node[above=1ex] {\ytext};

\foreach \z/\ztext in {1.25/$X_{\cF}^{++}$,3.75/$X_{\cF}^+$,6.25/$X_{\cF}$,8.75/$X_{\cE}$}
  \draw (\z,0) node[below] {\ztext};
\end{tikzpicture}
\caption{The slice of the movable cone $\cMov(X_\cF)$ for $M= V_5$.} \label{movV5sliPic}
\end{figure}

\begin{remark}
It is immediate that $X_\cE$ and $X_{\cF}^{++}$ are not isomorphic. Indeed, if they were isomorphic, then $Z_\cF$ and $D_1 (\sigma)$ would be isomorphic. But this is impossible because the number of singularities of $Z_\cF$ and $D_1 (\sigma)$ are $1$ and $29$ by Theorem \ref{mainV5} \eqref{mainV53} and Table \ref{intAll} respectively.
\end{remark}






\subsection{\texorpdfstring{Another description of $\bP (\cF)$}{Another description of P(F)}}\label{dP5_subsection1}


Let $W$ be a five-dimensional vector space over $\bC$.
Fix a two-dimensional subspace $ U \subseteq \wedge^3 W=H^0(\Gr(2,W),\cO(1))$ such that 
\[
    V_5=\Gr(2,W) \cap \bP(\wedge^3 W / U) \hookrightarrow  \bP(\wedge^3 W) 
\]
is a smooth del Pezzo fourfold of degree $5$.

Set $\cT = T_{ \bP(W^{\vee})} (-1)$. We define a coherent sheaf $\cC$ on $\bP(W^{\vee})$ as the cokernel of
\begin{equation}\label{eq_C'2}
\tau \colon U \otimes \cO_{ \bP(W^{\vee})} \hookrightarrow  \wedge^3 W \otimes \cO_{ \bP(W^{\vee})} \ra \wedge^3 \cT.
\end{equation}
Since $\bP_{\Gr(2,W)} ( \cS_{\Gr(2,W) }^{\vee} )$ is the flag variety $F(1,2 ; W )$, we have the following diagram: 
\begin{equation*}\label{V5twoprojs}
    \begin{tikzcd}
        & \bP_{\Gr(2,W)} ( \cS_{\Gr(2,W) }^{\vee} ) \ar[rd, "p_2"] \ar[ld, "p_1", swap]  &  \\
        \Gr(1,W) =\bP(W^{\vee}) && \Gr(2,W).      
    \end{tikzcd}
\end{equation*}
Here $p_1$ and $p_2$ are the two projection maps. We note that 
\begin{align}
    \bP_{\Gr(2,W)} ( \cS_{\Gr(2,W) }^{\vee} ) &= \bP_{\bP(W^{\vee}) }( \cT^{\vee}) \cong \bP_{\bP(W^{\vee}) }( \wedge^3 \cT ), \label{idfGrT}\\
    p_1^\ast \cO_{\bP (W^\vee)} (1) &= \cO_{\cS_{\Gr(2,W) }^{\vee}} (1). \label{idfGrTO1}
\end{align}

Since $ V_5 \subseteq \Gr(2,W) $ is the zero locus of $U \otimes \cO_{\Gr(2,W)} \ra \cO_{\Gr(2,W)} (1)$, the projective bundle $\bP_{V_5} ( \cS_{V_5}^{\vee} ) \cong \bP (\cF)$ is the zero locus of
\begin{equation}\label{eq_zero_F}
U \otimes \cO \ra  p_2^* \cO_{\Gr(2,W)} (1) 
\end{equation}
in $\bP_{\Gr(2,W)} ( \cS_{\Gr(2,W) }^{\vee} )$.



\begin{lemma}\label{V5anoth}
The projective bundle $\bP_{V_5} ( \cS_{V_5}^{\vee} )$ coincides with 
\[
    \bP_{\bP(W^{\vee})} (\cC) \subseteq  \bP_{\bP(W^{\vee}) }( \wedge^3 \cT )
\]
under the identification \eqref{idfGrT}.
\end{lemma}

\begin{proof}
We apply 
Lemma \ref{detzer} \eqref{detzer1} to $\tau \colon U \otimes \cO_{\bP(W^{\vee})} \to  \wedge^3 \cT  $ on $\bP(W^{\vee}) $.
Under the identification \eqref{idfGrT},
we have $p_2^* \cO_{\Gr(2,W)} (1) = \cO_{\wedge^3 \cT} (1) $.
Then 
\eqref{eq_zero_F} is nothing but the section $s_{\tau} $ induced by $\tau$,
and hence this lemma follows.
\end{proof}

\begin{remark}\label{rem_Banica}
By Lemma \ref{V5anoth}, the fiber of $\bP_{V_5} ( \cS_{V_5}^{\vee} ) = \bP_{\bP(W^{\vee})} (\cC) \to \bP(W^{\vee})$ over $x \in D_k(\tau) \setminus D_{k-1}(\tau)$ is $\bP^{4-k-1}=\bP^{3-k}$ for $k=0,1,2$.
Then we see that $ D_0(\tau)=\varnothing$ and $\dim D_1(\tau) =1$  by \cite[Proposition 1.8, Theorem 4.7]{Kan19}.
\end{remark}

\subsection{\texorpdfstring{Another description of $X_{\cF}$ and Construction of the flop $X^+_{\cF}$}{Another description of XF and Construction of the flop X+F}} \label{V5anflop_subsec}

We are going to give another description of $X_{\cF}$ using Lemma \ref{V5anoth}. By this description, we will construct a small contraction $X_{\cF} \to  Y_{\cF}$ and its flop $X^+_{\cF} \to  Y_{\cF}$.


Recall that, by Proposition \ref{nodal3F}, the $X_\cF \subseteq \bP_{V_5} (\cF) \cong \bP_{V_5}(\cS_{V_5}^{\vee})$ is induced by a general morphism $\sigma \colon \cE^{\vee} \to \cF$, i.e., $X_\cF = Z (s_\sigma)$. Observe that $\cE = \cO^{\oplus 2}$ and thus a general morphism $\sigma$ corrseponds to a general two-dimensional subspace of $H^0(\cF) =H^0 (\cS_{V_5}^{\vee} (1))$, say $U_{\sigma}$. Note that $H^i (\cS_{\Gr (2, W)}^{\vee} (1 - i)) = 0$ for $i = 1, 2$ and thus the natural map  
\[
    W^{\vee} \otimes \wedge^3 W = H^0 (W^{\vee} \otimes \cO_{\Gr (2, W)} (1)) \twoheadrightarrow H^0 (\cS_{\Gr (2, W)}^{\vee} (1)) \twoheadrightarrow H^0 (\cS_{V_5}^{\vee} (1))
\]
is surjective. 
Hence we can choose a two-dimensional subspace $U' \subseteq W^{\vee} \otimes \wedge^3 W$ which maps onto $U_{\sigma}$.
In other words,
we can take general $U' \subseteq W^{\vee} \otimes \wedge^3 W$ 
and identify $\sigma \colon \cE^\vee \to \cF$ with the natural morphism
\[
U' \otimes \cO_{V_5} \to W^{\vee} \otimes \cO_{V_5} (1) \to \cS_{V_5}^{\vee} (1) =\cF.
\]

Let $U ' \subseteq W^{\vee} \otimes \wedge^3 W$ be a general two-dimensional subspace as above. It induces $U' \otimes \cO \ra \cO_{\Gr(2,W)} (1)\boxtimes \cO_{\bP(W^{\vee})} (1) $ on $\bP_{\Gr(2,W)}(\cS_{\Gr(2,W)}^{\vee}) $
and $U' \otimes  \cO_{\bP(W^{\vee})}  \ra   \wedge^3 \cT (1)$ on $\bP(W^{\vee})$. We denote by
\begin{equation*}
    \tau' \colon U' \otimes  \cO_{\bP(W^{\vee})} (- 1)  \ra   \wedge^3 \cT
\end{equation*}
the twisted morphism. Set
\[
    \cH \coloneqq (U^{\vee} \otimes \cO_{ \bP(W^{\vee})}) \oplus  (U'^{\vee} \otimes  \cO_{\bP(W^{\vee})}(1) ),
\]
where $U \subseteq \wedge^3 W$ is the same two-dimensional subspace as in the beginning of Section \ref{dP5_subsection1}.



\begin{proposition}\label{V5CY3}
Let $\varrho \colon  \cH^{\vee}   \ra \wedge^3 \cT$ be the sum of the morphisms $\tau$ and $\tau'$.
\begin{enumerate}[(i)]
    \item The Calabi--Yau threefold $X_{\cF}  \subseteq \bP_{V_5}(\cS_{V_5}^{\vee})$ coincides with 
        \[
            \bP_{\bP(W^{\vee})} (\coker \varrho) \subseteq  \bP_{\bP(W^{\vee})} (\cC)
        \]
    under the identification \eqref{idfGrT} if $U'$ is chosen in general.\label{V5CY3_item1}
    \item The image $p_1(X_{\cF})$ coincides with $Y_{\cF} \coloneqq D_3(\varrho) $ and $\pi_{\wedge^3 \cT} \colon X_{\cF} \to Y_{\cF} \subseteq \bP (W^{\vee}) $ is a small resolution of a quintic hypersurface with $54$ ODPs. \label{V5CY3_item2}
\end{enumerate}
\end{proposition}

\begin{proof}
\eqref{V5CY3_item1}
Applying Lemma \ref{detzer} \eqref{detzer1} to $\varrho$ on $\bP(W^{\vee})$,
we have $  \bP_{\bP(W^{\vee})} (\coker \varrho) =Z(s_{\varrho}) $.
Hence it suffices to show $Z(s_{\varrho}) =X_{\cF}$.

To shorten notation, we write $\cU^{\vee}_Z$ instead of $U^{\vee} \otimes \cO_Z$ for a variety $Z$, and similarly for $\cU_Z'^{\vee}$. Then the section $s_{\varrho} \in H^0( \cH \boxtimes  \cO_{\wedge^3 \cT} (1))$ is the sum of 
\[
    s_{\tau} \in H^0( \cU_{\bP(W^{\vee})}^{\vee} \boxtimes \cO_{\wedge^3 \cT} (1) ) \quad \text{and} \quad
    s_{\tau'} \in H^0(\cU_{\bP(W^{\vee})}'^{\vee} (1) \boxtimes \cO_{\wedge^3 \cT} (1)).
\]
Under the identification 
$ \bP_{\Gr(2,W)} ( \cS_{\Gr(2,W) }^{\vee} ) = \bP_{\bP(W^{\vee}) }( \wedge^3 \cT ) $ in \eqref{idfGrT},
it holds that
\[
\cO_{\bP(W^{\vee})}(1) \boxtimes \cO_{\wedge^3 \cT} (1)
=p_1^\ast \cO_{\bP (W^\vee)} (1)\otimes \cO_{\wedge^3 \cT} (1) 
=\cO_{\cS_{\Gr(2,W) }^{\vee}} (1) \otimes p_2^* \cO_{\Gr(2,W)}(1) .
\]
Hence $s_{\tau'}$ is a section in $H^0( \cU_{\Gr(2,W)}'^{\vee} (1) \boxtimes \cO_{\cS_{\Gr(2,W) }^{\vee}} (1) )$ on $ \bP_{\Gr(2,W)} ( \cS_{\Gr(2,W) }^{\vee} )$.


By Lemmas \ref{detzer} and \ref{V5anoth},
we have $Z(s_{\tau}) =  \bP_{\bP(W^{\vee})} (\cC) =\bP_{V_5}(\cS_{V_5}^{\vee})$.
By the choice of $U'$,
the restriction of $s_{\tau'}$ on $\bP_{V_5}(\cS_{V_5}^{\vee})$, which is a section in 
\[
    H^0( \cU_{V_5}'^{\vee} (1) \boxtimes \cO_{\cS_{V_5 }^{\vee}} (1) ) 
    = H^0( \cU_{V_5}'^{\vee} \boxtimes \cO_{\cF} (1) ) ,
\]
can be identified with $s_{\sigma} \in H^0(\cE^{\vee} \boxtimes \cO_{\cF}(1))$ on $\bP_{V_5} (\cF) \cong \bP_{V_5}(\cS_{V_5}^{\vee})$.
Thus we have
\[
Z(s_{\varrho}) = Z(s_{\tau}) \cap Z(s_{\tau'})= \bP_{V_5}(\cS_{V_5}^{\vee}) \cap  Z(s_{\tau'}) = Z(s_{\tau'}|_{\bP_{V_5}(\cS_{V_5}^{\vee})}) = Z(s_{\sigma})=X_{\cF}
\]
and \eqref{V5CY3_item1} follows.\\
\eqref{V5CY3_item2}
By \eqref{V5CY3_item1} and Lemma \ref{detzer} \eqref{detzer2}, the restriction of  $p_1$ to $X_{\cF}$ maps onto $Y_{\cF}=D_3(\varrho)$ and the fiber of $\pi_{\wedge^3 \cT} \colon X_{\cF} \to Y_{\cF}$ over $x \in D_k(\varrho)\setminus D_{k-1}(\varrho) $ is $\bP^{4-k-1}=\bP^{3-k}$.
Thus $Y_\cF \setminus D_2 (\varrho) \cong X_\cF \setminus \pi_{\wedge^3 \cT}^{-1} (D_2 (\varrho))$ is smooth of dimension three and the smallness of $X_{\cF} \to Y_{\cF}$ holds if $D_1(\varrho)=\varnothing$ and $\dim D_2(\varrho) = 0$. Indeed, we shall prove that $\varrho$ is $3$-general. 

Recall that $D_k(\varrho)$ is the locus where the rank of $\coker \varrho$ is at least $4-k$.
Furthermore,
$\coker \varrho$ coincides with the cokernel of the composite
\[
    \tau'' \colon   U' \otimes  \cO_{\bP(W^{\vee})} (- 1)   \xrightarrow{\tau'}  \wedge^3 \cT \to \cC .
\]
Though $\cC = \coker \tau$ is not locally free, the restriction of $\cC$ on $D_{l}(\tau) \setminus D_{l-1}(\tau)$ is locally free of rank $4-l$.
Since the composite $(W^{\vee} \otimes \wedge^3 W) \otimes  \cO_{\bP(W^{\vee})}  \ra   \wedge^3 \cT (1)  \to \cC (1) $ is surjective,
so is the restriction of it on $D_{l}(\tau) \setminus D_{l-1}(\tau)$.
Hence we can apply Bertini-type theorem (cf.\ \cite[Theorem 3.1]{LW2022}) to compute the codimension of $D_{k}( \tau''|_{D_{l}(\tau) \setminus D_{l-1}(\tau)})$ in $D_{l}(\tau) \setminus D_{l-1}(\tau)$.

Note that $D_{2}(\tau)=\bP(W^{\vee}) $ and $\dim D_1(\tau)=1, D_0(\tau)=\varnothing$ by Remark \ref{rem_Banica}.
On $ D_{2}(\tau) \setminus D_{1}(\tau) = \bP(W^{\vee}) \setminus D_1(\tau)$,
the $\cC$ is locally free of rank two and hence $D_k(\varrho)  \setminus D_1(\tau)$ coincides with $D_{k-2}( \tau''|_{\bP(W^{\vee}) \setminus D_1(\tau)})$.
Thus $D_1(\varrho)  \setminus D_1(\tau) = \varnothing$
and $D_2(\varrho)  \setminus D_1(\tau) $ is smooth with at most zero-dimension by Bertini-type theorem.
On the other hand, $\cC$ is locally free of rank three on $ D_{1}(\tau) \setminus D_{0}(\tau)= D_1(\tau) $.
Hence $D_2(\varrho) \cap D_1(\tau)  $ coincides with $D_{1}( \tau''|_{D_1(\tau)})$.
Then the expected codimension of $D_{1}(\tau''|_{D_1(\tau) })  $ is two and hence 
$D_2(\varrho) \cap D_1(\tau)  =D_{1}(\tau''|_{D_1(\tau) })  = \varnothing$ by $\dim D_1(\tau)  =1$.

Thus we have $D_1(\varrho)=\varnothing$ and $D_2(\varrho) = D_2(\varrho)  \setminus D_1(\tau) $ is smooth with $\dim D_2(\varrho) \leqslant 0$. Note that $Y_\cF \subseteq \bP(W^{\vee})$ is a quintic hypersurface since $Y_{\cF}=D_3(\varrho)$ is the zero locus of $\det (\varrho)$, which is a global section of $\det (\cH) \otimes \det (\wedge^3 \cT) \cong \cO_{\bP (W^{\vee})} (5)$.
Then $\rho(Y_{\cF}) =1$ and  $\pi_{\wedge^3 \cT}$ is not an isomorphism by $\rho(X_{\cF}) =2$, i.e, $D_2(\varrho) \neq \varnothing$. Hence $\varrho$ is $3$-general and $\pi_{\wedge^3 \cT}$ is a small birational morphism. Applying Proposition \ref{LHinterODP}, we find that the number of ODPs of $Y_\cF$ is 
\begin{equation}\label{V5YFodp}
    |D_2 (\varrho)| = \int_{\bP (W^{\vee})} c_2 (\wedge^3 \cT - \cH^{\vee})^2 - c_1 (\wedge^3 \cT - \cH^{\vee}) \cdot c_3 (\wedge^3 \cT - \cH^{\vee}) = 54.
\end{equation}
This completes the proof.
\end{proof}



To find the flop of $X_{\cF}$, we can use the same trick as in Section \ref{Nefcone_subsec}. Let
\begin{align*}
    X^+_{\cF} \coloneqq \bP_{\bP(W^{\vee})}(\coker \varrho^{\vee}  )  \subseteq  \bP_{\bP(W^{\vee})}(\cH) ,
\end{align*}
Note that $D_k (\varrho) =  D_k (\varrho^\vee)$ for all $k$. Therefore $\varrho^\vee$ is also $3$-general. By Proposition \ref{n-general}, we find that the determinantal contraction $\pi_{\cH} \colon X^+_{\cF}  \ra Y_{\cF}$ is small and the $X^+_{\cF}$ is smooth. Let $H'$ be the pullback of $\cO_{\bP(W^{\vee}) } (1) $ via $X^+_{\cF}  \ra \bP(W^{\vee})$ and $L' =\cO_{\cH}(1) |_{X^+_{\cF} } $. Applying Proposition \ref{LHinterODP}, we get Table \ref{intV5‘}. 

\begin{table}
  \centering
  \begin{tabular}{ccccccc}
    \toprule
    $L'^3$ & $L'^2.H'$ & $L'.H'^2$ & $H'^3$ & $L'.c_2(T_{X_{\cF}^+})$ & $H'.c_2(T_{X_{\cF}^+})$ & $\#$ of ODPs \\
    \midrule
     34 & 23 & 13 & 5 & 76 & 50 & 54 \\
    \bottomrule
  \end{tabular}
  \caption{The intersection numbers of $X_\cF^+$.}
  \label{intV5‘}
\end{table}




\subsection{Construction of another birational model \texorpdfstring{$X^{++}_{\cF}$}{X++F}}

In the previous subsection,
we construct a small contraction $\pi_{\cH} : X^+_{\cF} \to Y_{\cF} \subseteq \bP(W^{\vee})$.
We now turn to find
another small contraction $X^{+}_{\cF} \to Z_{\cF}$ and its flop $X^{++}_{\cF}$.

Recall that $ U \subseteq \wedge^3 W$ is a two-dimensional subspace such that 
$V_5=\Gr(2,W) \cap \bP(\wedge^3 W / U) $ is the smooth del Pezzo fourfold $V_5$ of degree $5$.
By $\wedge^3 W \cong \wedge^2 W^{\vee} $,
we might think that $U \subseteq \wedge^2 W^{\vee}$.

Each $\omega \in  \wedge^2 W^{\vee}  $ defines an alternating form $ f_{\omega } \colon W \times W \ra \bC$
with the kernel
\[
\ker f_{\omega} =\{ w \in W \, \mid \, f_{\omega}(w,w') =0 \text{ for any } w' \in W \}.
\]
Since $ \rank f_{\omega}$ is even and $\dim W=5$,
the dimension of $\ker f_{\omega} $ is $1$ or $3$ if $\omega \neq 0$.

\begin{lemma}\label{lem_ker_f_u}
For any non-zero $ u \in U$, we have $\dim \ker f_u =1$.
\end{lemma}

\begin{proof}
For $\omega \in  \wedge^2 W^{\vee} $,
let $H_{\omega} \subseteq \bP( \wedge^2 W^{\vee} )=\bP(\wedge^3 W) $
be the hyperplane defined by $\omega$.
Let $x= [V] \in \Gr(2,W)$ be a point corresponding to a two-dimensional subspace $V \subseteq W$.
Then it is known that $ H_{\omega}$ is tangent to $\Gr(2,W) \subseteq \bP( \wedge^2 W^{\vee} )$ at $x$
if and only $ V \subseteq \ker f_{\omega}$ (see \cite[Proposition 1.5]{BC09} for example).
Hence the singular locus of $ \Gr(2,W) \cap H_{\omega}$ coincides with
\[
\{ [V] \in \Gr(2,W) \, | \, V \subseteq \ker f_{\omega} \},
\]
which is empty or isomorphic to $\bP^2$ according to  $\dim \ker f_{\omega}=1$ or $ 3$.

If $\dim \ker f_u =3$ for some non-zero $ u \in U$, 
the $ \Gr(2,W) \cap H_{u}$ is singular along $\bP^2$.
Since $V_5$ is a linear section of $ \Gr(2,W) \cap H_{u}$ of codimension one, 
the $V_5$ is singular along $\bP^1$, which contradicts to the smoothness of $V_5$.
Hence we have $\dim \ker f_u =1$ for any non-zero $ u \in U$.
\end{proof}

By $\wedge^3 W \cong \wedge^2 W^{\vee} $ and $\wedge^3\cT \cong  \cT^{\vee}(1) $,
we might think that $\cC$ is the cokernel of 
\begin{equation} \label{eq_C'3}
\tau \colon U \otimes \cO_{ \bP(W^{\vee})} \hookrightarrow  \wedge^2 W^{\vee}  \otimes \cO_{ \bP(W^{\vee})} \ra \cT^{\vee}(1) .
\end{equation}

Set $\Sigma \coloneqq D_1 ( \tau) = D_1 ( \tau^{\vee}) \subseteq \bP(W^{\vee})$.
By \eqref{Cayleytrick}, we get the zero scheme
\[
    \Sigma^+ =\bP_{\bP(W^{\vee})}(\coker \tau^{\vee}) \subseteq \bP_{\bP(W^{\vee})}(U^{\vee} \otimes \cO_{\bP(W^{\vee})}) = \bP(U^{\vee} ) \times \bP(W^{\vee})
\] 
induced by $\tau^{\vee}$ and a natural morphism $\Sigma^+ \to \Sigma$.
We can also regard $\Sigma^+ \subseteq \bP_{\bP(W^{\vee})}(\cH) $
by the natural inclusion  $ \bP_{\bP(W^{\vee})}(U^{\vee} \otimes \cO_{\bP(W^{\vee})})  \subseteq  \bP_{\bP(W^{\vee})}(\cH)$.
By Remark \ref{rem_Banica}, we know that $\dim \Sigma=1$.
In fact, we have an explicit description of $\Sigma$ as follows:

\begin{lemma}\label{lem_isom_over_Z}
With notation as above:
\begin{enumerate}[(i)]
\item\label{lem_isom_over_Z1} The natural projection $\Sigma^+ \to \bP(U^{\vee}) \cong \bP^1$ is an isomorphism.
\item\label{lem_isom_over_Z2} The natural projection $\Sigma^+ \to  \Sigma$ is an isomorphism and it holds that
\begin{equation} \label{p1cur}
\Sigma = \{ [(\ker f_u)^{\vee}] \in \bP(W^{\vee}) \mid u \in U \setminus \{0\} \}.
\end{equation}
Furthermore, the degree of $\Sigma$ in $\bP (W^\vee)$ is two. 
\item\label{lem_isom_over_Z3} The $\Sigma^+ =  X^+_{\cF} \cap  \bP_{\bP(W^{\vee})}(U^{\vee} \otimes \cO_{\bP(W^{\vee})})  $ holds scheme-theoretically,
where we take the intersection in $\bP_{\bP(W^{\vee})}(\cH)  $.
\end{enumerate}
\end{lemma}

\begin{proof}
\eqref{lem_isom_over_Z1} 
For simplicity,
we sometimes denote $W \otimes \cO_Y $ by $\mathcal{W}_Y$ for a variety $Y$.
The subspace $U \subseteq \wedge^2 W^{\vee} $ induces a two-form
$  \wedge^2 \mathcal{W}_{ \bP(U^{\vee} ) } \to  \cO_{\bP(U^{\vee})}(1)$ on $\bP(U^{\vee})$.
Let $\mathcal{K} \subseteq \mathcal{W}_{\bP(U^{\vee})}$ be the kernel of this two-from,
that is, the kernel of the induced morphism
$ \mathcal{W}_{\bP(U^{\vee})} \to \mathcal{W}_{\bP(U^{\vee})}^{\vee} \otimes  \cO_{\bP(U^{\vee})}(1)$.
By Lemma \ref{lem_ker_f_u}, 
the $\mathcal{K} \subseteq \mathcal{W}_{\bP(U^{\vee})}$ is a subbundle  of rank one and
$\mathcal{K} ([u^{\vee}]) = \ker f_u \subseteq W$ holds for any non-zero $u \in U$.
Here $[u^{\vee}] \in \bP(U^{\vee})$ is the point corresponding to the quotient $ U^{\vee} \twoheadrightarrow (\bC u)^{\vee}$.

Let $p_{\bP(U^{\vee} )}$ and $ p_{\bP(W^{\vee})}$ be the natural projections from $\bP(U^{\vee} ) \times \bP(W^{\vee}) $ to the first and second factors respectively.
By Lemma \ref{detzer} \eqref{detzer1}, 
the $\Sigma^+ \subseteq \bP(U^{\vee} ) \times \bP(W^{\vee})$ is the zero locus of 
\begin{equation*} 
 s_{\tau^{\vee} } \colon  p_{\bP(W^{\vee})}^* (\cT(-1) ) \to  \wedge^2 \mathcal{W}_{ \bP(U^{\vee} ) \times \bP(W^{\vee})} \to
 p_{\bP(U^{\vee} )}^*  \cO_{\bP(U^{\vee})}(1) .
\end{equation*}
By the natural exact sequence
\begin{equation} \label{eq_wedgEu}
    0 \ra   \cT(-1) \ra \wedge^2 W  \otimes \cO_{ \bP(W^{\vee})} \ra  \wedge^2 \cT \ra 0
\end{equation}
obtained from $0 \ra \cO_{ \bP(W^{\vee})} (-1) \ra W \otimes \cO_{ \bP(W^{\vee})} \ra  \cT \ra 0$,
the zero locus of $ s_{\tau^{\vee} } $ is the maximum closed subscheme on which 
the two-form
$ \wedge^2 \mathcal{W}_{ \bP(U^{\vee} ) \times \bP(W^{\vee})} \to p_{\bP(U^{\vee} )}^*  \cO_{\bP(U^{\vee})}(1)$
factors through $\wedge^2 \mathcal{W}_{ \bP(U^{\vee} ) \times \bP(W^{\vee})} \to \wedge^2 (p_{\bP(W^{\vee})}^* \cT ) $.
By the definition of the kernel $\mathcal{K}$ and $\cT=  \mathcal{W}_{ \bP(W^{\vee})}/\cO_{ \bP(W^{\vee})} (-1) $,
the zero locus $ \Sigma^+ =Z(s_{\tau^{\vee} } )$ is nothing but 
the maximum locus on which $ p_{\bP(W^{\vee})}^* \cO_{ \bP(W^{\vee})} (-1) \subseteq  p_{\bP(U^{\vee})}^*\mathcal{K}$,
that is, 
the zero locus of 
\[
p_{\bP(W^{\vee})}^* \cO_{ \bP(W^{\vee})} (-1)  \hookrightarrow  \mathcal{W}_{ \bP(U^{\vee} ) \times \bP(W^{\vee}) } 
\to p_{\bP(U^{\vee})}^* (\mathcal{W}_{\bP(U^{\vee})} /\mathcal{K}).
\]
By Lemma \ref{detzer} \eqref{detzer1},
this zero locus coincides with $\bP_{\bP(U^{\vee})}(\mathcal{K}^{\vee})$.
Since $\mathcal{K}$ is locally free of rank one,
the projection $\Sigma^+=\bP_{\bP(U^{\vee})}(\mathcal{K}^{\vee}) \to \bP(U^{\vee})$ is an isomorphism.

\noindent
\eqref{lem_isom_over_Z2} 
By the proof of \eqref{lem_isom_over_Z1}, we know that
\begin{align}\label{lem_subK}
    p_{\bP(W^{\vee})}^* \cO_{ \bP(W^{\vee})} (-1)|_{\Sigma^+} \subseteq  p_{\bP(U^{\vee})}^*\mathcal{K}|_{\Sigma^+} 
\end{align}
holds on $\Sigma^+$. Since both sides of \eqref{lem_subK} are subbundles of $ \mathcal{W}_{ \bP(U^{\vee} ) \times \bP(W^{\vee})}$ of the same rank, this shows that \eqref{lem_subK} is actually an equality.
Thus 
$\bP(U^{\vee}) \cong \Sigma^+ \to \Sigma \subseteq \bP(W^{\vee})$
is induced by the quotient $ \mathcal{W}_{\bP(U^{\vee})}^{\vee} =W^{\vee} \otimes \cO_{\bP(U^{\vee})}  \to \mathcal{K}^{\vee} $
and hence $y=[u^{\vee}] \in \bP(U^{\vee}) $ is mapped to $[\mathcal{K}^{\vee}(y)] = [(\ker f_u)^{\vee}] \in \Sigma $.
Since $\Sigma^+ \to \Sigma=D_1(\tau^{\vee}) $ is surjective by Lemma \ref{detzer} \eqref{detzer2},
we obtain \eqref{p1cur}. 


By Lemma \ref{detzer} \eqref{detzer2}, the map $\Sigma^+ \to \Sigma$ is an isomorphism if $D_0(\tau^{\vee})=\varnothing$. If $D_0(\tau^{\vee})$ is non-empty, then the fiber of $\Sigma^+= \bP_{\bP(W^{\vee})}(\coker \tau^{\vee})  \to \Sigma$ over a point $x \in D_0(\tau^{\vee}) $ is $\bP^1$,
that is, the fiber coincides with the whole $\Sigma^+ \cong \bP(U^{\vee})=\bP^1$.
Therefore the surjective morphism $\Sigma^+ \to \Sigma$ is an isomorphism if it is not a constant map.
We already saw that $\bP(U^{\vee}) \cong \Sigma^+ \to \Sigma \subseteq \bP(W^{\vee})$ is induced by
$ W^{\vee} \otimes \cO_{\bP(U^{\vee})} \twoheadrightarrow \mathcal{K}^{\vee} $.
Hence $\Sigma^+ \to \Sigma$ is not a constant map if and only if $ \mathcal{K}^{\vee}  \not \cong \cO_{\bP(U^{\vee})}$.
In fact, we can show $\mathcal{K}^{\vee} \cong \cO_{\bP(U^{\vee})}(2) $ as follows:
By the definition of $\mathcal{K}$,
the two-form
$  \wedge^2 \mathcal{W}_{ \bP(U^{\vee} ) } \to  \cO_{\bP(U^{\vee})}(1)$ induces
\begin{equation} \label{eq_nondegenerate_two_form}
  \wedge^2 (\mathcal{W}_{ \bP(U^{\vee} ) }/\mathcal{K}) \to  \cO_{\bP(U^{\vee})}(1),
\end{equation}
which is non-degenerate at any point in $\bP(U^{\vee})$.
Hence it induces an isomorphism $ \mathcal{W}_{ \bP(U^{\vee} ) }/\mathcal{K} \to (\mathcal{W}_{ \bP(U^{\vee} ) }/\mathcal{K})^{\vee} \otimes  \cO_{\bP(U^{\vee})}(1) $.
By taking the determinant,
we see that $\mathcal{K}^{\vee} \cong \cO_{\bP(U^{\vee})}(2) $.

Thus $\Sigma^+ \to \Sigma$ is an isomorphism induced by $\mathcal{K}^{\vee} \cong \cO_{\bP(U^{\vee})}(2) $.
In particular, the degree of $\Sigma $ is two\footnote{There is an alternative way to compute the degree of $\Sigma = D_1 (\tau) \subseteq \bP (W^\vee) $: 
\[
    \Sigma \cdot c_1 (\cO_{\bP (W^\vee)} (1)) = c_3 (\wedge^3 \cT) \cdot c_1 (\cO_{\bP (W^\vee)} (1)) = 2
\]
by \eqref{eq_C'2} and \cite[Theorem 14.4 (c)]{fulton98}.}

\noindent
\eqref{lem_isom_over_Z3} 
Recall $\varrho = (\tau, \tau')$.
By definition and Lemma \ref{detzer} \eqref{detzer1},
we have $\Sigma^+ = \bP_{\bP(W^{\vee})} (\coker \tau^{\vee}) =Z(s_{\tau^{\vee}}) $ and 
$X^+_{\cF} = \bP_{\bP(W^{\vee})} (\coker \varrho^{\vee}) =Z( s_{\varrho^{\vee}})$.
Since the restriction of $ s_{\varrho^{\vee}}=( s_{\tau^{\vee}}, s_{\tau'^{\vee}})$ on $ \bP_{\bP(W^{\vee})}(U^{\vee} \otimes \cO_{\bP(W^{\vee})}) $ 
is $ (s_{\tau^{\vee}} ,0)$,
we have 
\[
X^+_{\cF} \cap \bP_{\bP(W^{\vee})}(U^{\vee} \otimes \cO_{\bP(W^{\vee})}) = Z( s_{\varrho^{\vee}}) \cap \bP_{\bP(W^{\vee})}(U^{\vee} \otimes \cO_{\bP(W^{\vee})}) = Z(s_{\tau^{\vee}}) =\Sigma^+,
\]
which completes the proof.
\end{proof}

We will construct a small contraction $X^{+}_{\cF} \to Z_{\cF}$ whose exceptional locus is the curve $\Sigma^+$.
The task is now to find a divisor $D$ of $X_{\cF}^+$ not meeting $\Sigma^{+}$. Using the description \eqref{eq_C'3} of $\tau$ and similarly for $\tau'$, we find that $\varrho$ is decomposed as
\[
    \varrho \colon  \cH^{\vee} \ra \wedge^2 W^{\vee}  \otimes \cO_{ \bP(W^{\vee})} \ra \cT^{\vee} (1).
\]
On $p_{\cH} \colon  \bP_{\bP(W^{\vee})}(\cH) \ra \bP(W^{\vee})$, we have
\[
    \cO_{\cH} (-1) \hookrightarrow p_{\cH}^* \cH^{\vee} \ra \wedge^2 W^{\vee} \otimes \cO \ra p_{\cH}^* \cT^{\vee} (1),
\]
whose zero locus is nothing but $X^+_{\cF} =\bP_{\bP(W^{\vee})}(\coker \varrho^{\vee}  )  \subseteq \bP_{\bP(W^{\vee})}(\cH) $.
Hence the restriction of $ \cO_{\cH} (-1) \hookrightarrow p_{\cH}^* \cH^{\vee} \ra \wedge^2 W^{\vee} \otimes \cO $ on 
$X^+_{\cF}$ factors as 
\[
    s \colon \cO_{\cH} (-1)|_{ X^+_{\cF} }   \hookrightarrow  \pi_{\cH}^* (\cH^{\vee}) \ra \pi_{\cH}^* (\wedge^2 \cT^{\vee})
    \subseteq \wedge^2 W^{\vee} \otimes \cO
\]
by the dual of \eqref{eq_wedgEu},
where we recall that $\pi_{\cH} \colon X^+_{\cF}  \ra Y_{\cF} \subseteq \bP(W^{\vee})$ is the restriction of $p_{\cH}$.
By taking wedge product of $s$, we have
\[
    s \wedge s \colon (\cO_{\cH} (-1)  \otimes \cO_{\cH} (-1) )|_{ X^+_{\cF} }  \ra  \pi_{\cH}^*( \wedge^2 \cT^{\vee} \otimes  \wedge^2 \cT^{\vee}) \ra  \pi_{\cH}^* (\wedge^4 \cT^{\vee}).
\]
From the isomorphism $\wedge^4 \cT^{\vee} \cong \cO_{\bP(W^{\vee})}(-1)$, we get a section 
\[
    s \wedge s \in H^0( X^+_{\cF} , \cO_{\cH} (2)|_{ X^+_{\cF} }  \otimes \pi_{\cH}^* \cO_{\bP(W^{\vee})}(-1)).
\]
Let $D \subseteq X^+_{\cF} $ denote the zero locus of the section $s \wedge s$.
By the following lemma,
the $s \wedge s  $ is a non-trivial section and hence $D$ is a true (effective) divisor in the linear system $\left|\cO_{\cH} (2)|_{ X^+_{\cF} }  \otimes \pi_{\cH}^* \cO_{\bP(W^{\vee})}(-1) \right|$.

\begin{lemma}\label{lem_Z capD}
Let $\Sigma^+$ be as in Lemma \ref{lem_isom_over_Z}. Then $D \cap \Sigma^{+} = \varnothing$. 
\end{lemma}

\begin{proof}
Notice that the zero locus of $s \wedge s$ is the locus where the corresponding two-form $ s  \colon
\cO_{\cH} (-1)|_{ X^+_{\cF} } \to \pi_{\cH}^* (\wedge^2 \cT^{\vee})$
is degenerate. 

We use the notation in the proof of Lemma \ref{lem_isom_over_Z}.
Recall that an isomorphism $\bP(U^{\vee}) \to \Sigma$ is defined by $[u^{\vee}] \mapsto [(\ker f_u)^{\vee}]$
and $\Sigma^+ \subseteq \bP(U^{\vee}) \times \bP(W^{\vee})$ coincides with the graph of the isomorphism.
By the isomorphism $\bP(U^{\vee}) \to \Sigma$, the restriction $\cT|_{\Sigma}$ is identified with 
$\mathcal{W}_{\bP(U^{\vee})} / \mathcal{K} =(W \otimes \cO_{\bP(U^{\vee})}) / \mathcal{K} $.

Under the identification $\Sigma^+ \cong \bP(U^{\vee}) $,
the restriction of $s$ on $\Sigma^+$ is identified with
\[
    \cO_{\bP(U^{\vee}) } (-1)    \ra \wedge^2 (\mathcal{W}_{\bP(U^{\vee})} / \mathcal{K}),
\]
which is nothing but the dual of \eqref{eq_nondegenerate_two_form}.
It is non-degenerate at any point in $\bP(U^{\vee}) $
by the definition of the kernel $\mathcal{K}$.
Thus $s$ is non-degenerate at any point in $\Sigma^+$ and hence
$\Sigma^+ $ does not intersect with the zero locus of $s \wedge s$. 
\end{proof}

Recall that $L' =\cO_{\cH}(1) |_{X^+_{\cF} } $ and $H'$ is the pullback of $\cO_{\bP(W^{\vee}) } (1) $ via $X^+_{\cF}  \ra \bP(W^{\vee})$.
Using the effective divisor $D$, we can construct a small contraction $X^+_{\cF} \to Z_{\cF}$
and its flop $X^{++}_{\cF}$ as follows:

\begin{proposition}\label{V5XF+Nef}
There exists a birational model $X^{++}_{\cF} $ of $ X_{\cF}^+$ such that
\begin{align*}
\Nef (  X_{\cF}^+) &= \bR_{ \geqslant 0} [H']+ \bR_{ \geqslant 0}   [2L'-H'] ,\\
\Nef (  X_{\cF}^{++}) &= \bR_{ \geqslant 0}  [2L'-H'] +   \bR_{ \geqslant 0}  [L'-H'].
\end{align*}
Furthermore,
$L'-H' $ and $2 L' - H'$ define a fibration $ X^{++}_{\cF} \ra \bP^1$ and a small contration $\psi \colon X_{\cF}^+ \to Z_{\cF}$ with $\Exc (\psi) = \Sigma^+$ being contracted to an ODP respectively.
\end{proposition}

\begin{proof}
Note that
\[
    H^0\left(\cO_{\cH} (1)  \otimes p_{\cH}^*  \cO_{\bP(W^{\vee})}(-1) \right) = H^0 \left(\cH^{\vee} \otimes\cO_{\bP(W^{\vee})}(-1) \right)= U'^{\vee}.
\]
Take a basis $t_1,t_2 \in U'^{\vee}$ and 
let $D_1',D_2'  \subseteq \bP_{\bP(W^{\vee})}(\cH) $ be the corresponding divisors. 
Then  $D_1' \cap D_2' =  \bP_{\bP(W^{\vee})}(U^{\vee} \otimes \cO_{\bP(W^{\vee})}) \subseteq \bP_{\bP(W^{\vee})} (\cH) $.
Hence $X^+_{\cF} \cap D_1' \cap D_2' =\Sigma^+$ by Lemma \ref{lem_isom_over_Z}.
By Lemma \ref{lem_Z capD},
we have 
$X^+_{\cF} \cap D_1' \cap D_2' \cap D =\Sigma^+ \cap D =\varnothing$. Set $D_i \coloneqq D_i'|_{X^+_{\cF} }$ for $i = 1, 2$.
Hence, on $X^+_{\cF}$, we have effective divisors
\begin{center}
    $D_1, D_2 \in |L' - H'|$ \quad and \quad  $D \in |2 L' - H'|$.
\end{center}
Since $D_1\cap  D_2 \cap D =\varnothing$,
we can apply \cite[Lemma 2.4, Proposition 2.5]{Ito14} and 
obtain a model $X^{++}_{\cF} $.
The description of the nef cones and the morphism $ X^{++}_{\cF} \ra \bP^1$ is obtained  by \cite[Proposition 2.5]{Ito14}.

We review the construction of the flop for the benefit of the reader. First, since $H'$ is base point free, we have
\[
    \rBs (|2 L' - H'|) \subseteq D_1 \cap  D_2 \cap D =\varnothing.
\]
This means that $2 L' - H'$ is also base point free. Let $\psi \colon X_{\cF}^+ \to Z_{\cF}$ be the morphism defined by $|m (2 L' - H')|$ for $m \gg 0$. We claim that the exceptional set of $\psi$ is equal to $\Sigma^+$, i.e. $\psi$ is small.

Indeed, the curve $\Sigma^+$ is contained in $\Exc (\psi)$ since $\Sigma^+ \cdot (2 L' - H') = \Sigma^+ \cdot D = 0$ by Lemma \ref{lem_Z capD}. On other hand, if there exists a curve $C \subseteq \Exc (\psi)$ such that $C \neq \Sigma^+$, then $C \cdot (L' - H')$ and $C \cdot H'$ are nonnegative because $H'$ is base point free and $D_1 \cap D_2 = \Sigma^+ \neq C$. Therefore $C \cdot (2 L' - H') = 0$ implies $C \cdot L' = C \cdot H' = 0$. This contradicts the fact that $X_{\cF}^+$ is projective and has Picard number two.  

Next we construct the flop $X^{++}_{\cF}$ of $\psi \colon X_{\cF}^+ \to Z_{\cF}$.
Let $\mu \colon \Bl_{\Sigma^+} X_{\cF}^+ \to X_{\cF}^+$ be the blow-up along $\Sigma^+$ and $E$ the exceptional divisor of $\mu$. Since $D_1 \cap D_2 = \Sigma^+ \subseteq X_{\cF}^+$, the rational map $X_{\cF}^+ \dashrightarrow \bP (U'^{\vee})$ defined by the pencil $\left<D_1, D_2\right> \subseteq |L' - H'|$ can be resolved by $\mu$ to get the morphism $f \colon \Bl_{\Sigma^+} X_{\cF}^+ \to \bP (U'^{\vee}) \cong \bP^1$. Note that 
\[
    E = \bP_{\Sigma^+} (\cO (-D_1)|_{\Sigma^+} \oplus \cO (-D_2)|_{\Sigma^+})= \bP^1 \times \Sigma^+ \cong \bP^1 \times \bP^1.
\]
Furthermore, we have $\cO(-E)|_E \cong \cO_{\bP^1 \times \bP^1}(1,1)$
since $\Sigma^+ \cdot D_i  =-1$ by $\Sigma^+  \cdot (2L'-H')=0 $ and $\Sigma^+  \cdot H' =\Sigma  \cdot \cO_{\bP(W^{\vee})} (1) =2 $.
Hence $ \cO(\mu^* D_i -E) |_{E} \cong \cO_{\bP^1 \times \bP^1}(1,0)$
and
$f|_E \colon E=\bP^1 \times \Sigma^+ \to \bP^1 $ is the projection to the first factor.

We write $X_{\cF}^{++}$ for the normalization of the image of 
\[
    (f, \psi \circ \mu) \colon \Bl_{\Sigma^+} X_{\cF}^+  \to \bP^1 \times Z_{\cF},
\]
and the natural projection gives rise to the fibration $X_{\cF}^{++} \to \bP ^1$. 
Since $\psi \circ \mu(E) =\psi (\Sigma^+) \in Z_{\cF}$ is a point,
$(f, \psi \circ \mu)|_{E} $ contracts the divisor $E=\bP^1 \times \Sigma^+ $ to 
a curve $\bP^1 \times \psi(\Sigma^+)$. 
On the other hand, 
$ (f, \psi \circ \mu) $ is an isomorphim onto the image outside $E$
since $\Bl_{\Sigma^+} X_{\cF}^+ \setminus E \cong X_{\cF}^+ \setminus \Sigma^+ \cong Z_{\cF} \setminus \psi(\Sigma^+) $.
Therefore the natural projection $X_{\cF}^{++} \to Z_{\cF}$ is the flop of $X_{\cF}^{+} \ra Z_{\cF}$. In fact, it is an Atiyah flop because of $\cO (E)|_E \cong \cO_{\bP^1 \times \bP^1} (-1, -1)$.  

By construction,
$H'$ and $2L'-H'$ define small contractions $X^+_{\cF} \to Y_{\cF} $ and $X^+_{\cF} \to Z_{\cF} $ respectively.
Hence
$H'$ and $2L'-H'$ span $\Nef(X^+_{\cF})$.
On the other hand,
$2L'-H' $ and $ L' - H'$ define a small contration $ X_{\cF}^{++} \to Z_{\cF}$ and
a fibration $ X^{++}_{\cF} \ra \bP^1$ respectively.
Hence
$2L'-H' $ and $L' - H'$ span $\Nef(X^{++}_{\cF})$.
\end{proof}



\begin{proof}[Proof of Theorem \ref{mainV5}]
According to \eqref{idfGrTO1} and $\cS_{V_5}^{\vee} = \cF (-1)$, it follows that the determinantal contraction $\pi_{\wedge^3 \cT} \colon X_\cF \to Y_\cF$ is induced by $L - H$. Then \eqref{mainV52} follows from this and \eqref{V5YFodp}.

The \eqref{mainV51} and \eqref{mainV55} follow from the construction of $\pi_\cF$ and Proposition \ref{movXE}. Remark that, as in \eqref{K324}, we find that $c_2 (T_{V_5}) \cdot c_2 (\cF) = 53$ and $c_2 (T_{X_{\cE}}) \cdot L_\cE = 24$ by the Schubert calculus for $\Gr (2, W)$.

Set $\theta \coloneqq (\pi_{\wedge^3 \cT})^{- 1} \circ \pi_{\cH}$ and $\theta_\ast L' = a L + bH$. Notice that $\theta_\ast H' = L -H$, and $\theta$ is not an isomorphism because the numbers of singular points of $D_1 (\sigma)$, $Y_\cF$ and $Z_\cF$ are distinct. By Lemma \ref{ab}, Table \ref{intAll} and \ref{intV5‘}, we have
\[
    \begin{cases}
      23=47a^2+60ab+18b^2,\\
      13=17a+12b.\\
    \end{cases}
\]
The only solution is $(a, b) = (5, -6)$ because $ab < 0$. Then \eqref{mainV53} and \eqref{mainV54} follow from Proposition \ref{V5XF+Nef} except that a general fiber $S$ of $X_{\cF}^{++} \to \bP^1$ is K3. To prove this, we shall compute $c_2 (T_{X_{\cF}^{++}}) \cdot S$ as before. By Proposition \ref{V5XF+Nef}, \cite[Lemma 7.5]{Friedman91} and Table \ref{intV5‘}, the intersection number equals 
\[
    c_2 (T_{X_{\cF}^{+}}) \cdot (L' - H') + 2 (L' - H') \cdot \Sigma^+ = 26 - H' \cdot \Sigma^+.
\]
Notice that $H' \cdot \Sigma^+ =\cO_{\bP(W^{\vee})}(1) \cdot \Sigma =2 $
since $\Sigma \subseteq \bP (W^\vee)$ is a conic by Lemma \ref{lem_isom_over_Z}. Hence the intersection number $c_2 (T_{X_{\cF}^{++}}) \cdot S = 24$ and the proof is completed.
\end{proof}


\section{\texorpdfstring{The Grassmannian $\Gr (2, 4)$}{The Grassmannian Gr(2, 4)}}\label{G24_sec}

We consider the case $(M, \cF, \cE) = (\Gr(2, 4), \cS (2) \oplus \cO (1), \cO^{\oplus 3})$. As before, for a general morphism $\sigma \colon \cE^\vee \to \cF$, there are two smooth Calabi-Yau threefolds $X_{\cF}$, $X_{\cE}$ and the diagram \eqref{flopdiag}.
In this section,
we prove the following theorem.


\begin{theorem}\label{mainG24}
Let $(M, \cF, \cE) = (\Gr(2, 4), \cS (2) \oplus \cO (1), \cO^{\oplus 3})$. Then for a general morphism $\sigma \colon \cE^{\vee} \to \cF$, the scheme $X_{\cF}$ is a smooth Calabi--Yau threefold of Picard number two with
\[
    \Nef(X_\cF)=\bR_{\geqslant 0}[L-H]+\bR_{\geqslant 0}[H],
\]
such that
\begin{enumerate}[(i)]
    \item\label{mainG241} the determinantal contraction $\pi_{\cF}$ is induced by $H$;

    \item\label{mainG242} the $L -  H$ induces a small contraction $\varphi \colon X_{\cF} \to Y_{\cF}$ and $Y_{\cF} \subseteq \bP^4$ is a quintic hypersurface;
    
    \item\label{mainG243} the flop $X_\cF^{+}$ of $\varphi$ admits an elliptic fibration over $\bP^2$ induced by $4 L - 5 H$;

    \item\label{mainG244} the $X_{\cE}$ admits an elliptic fibration over $\bP^2$ induced by $4 H - L$.
\end{enumerate}
Moreover, the movable cone $\cMov(X_{\cF})$ is the convex cone generated by the divisors $4 L - 5 H$ and $4 H - L$ which is covered by the nef cones of $X_{\cF}$, $X_{\cF}^+$ and $X_{\cE}$, and there are no more minimal models of $X_{\cF}$, which we summarize in the following diagram:
\begin{equation*}
    \begin{tikzcd}
        X_{\cF}^+ \ar[d] \ar[rd, "\varphi^+"] \ar[rr, dashrightarrow, "\theta"] && X_{\cF} \ar[ld, "\varphi", swap] \ar[rd, "\pi_{\cF}"] \ar[rr, dashrightarrow, "\chi"] && X_{\cE} \ar[ld, "\pi_{\cE}", swap] \ar[d] \\
      \bP^2 & Y_{\cF} && D_2 (\sigma) & \bP^2
    \end{tikzcd}
\end{equation*}
\end{theorem}

The slice of the chamber structure of $\cMov(X_\cF)$ is given in Figure \ref{movG24sliPic}. We depict $X_{\cE}$, $X_{\cF}$, and $X_{\cF}^+$ inside their nef cones.

\begin{figure}[hbtp]
\centering
\begin{tikzpicture}
\draw(0,0)--(9,0);

\foreach \x in {0,3,6,9}
  \draw[very thick, teal] (\x,5pt)--(\x,-5pt);

\foreach \y/\ytext in {-0.25/$4L - 5H$,3/$L - H$,6/$H$,9/$4 H - L$}
  \draw (\y,0) node[above=1ex] {\ytext};

\foreach \z/\ztext in {1.5/$X_{\cF}^+$,4.6/$X_{\cF}$,7.5/$X_{\cE}$}
  \draw (\z,0) node[below] {\ztext};
\end{tikzpicture}
\caption{The slice of the movable cone $\cMov(X_\cF)$ for $M=\Gr (2, 4)$.} \label{movG24sliPic}
\end{figure}





\begin{remark}
It is easily seen that $X_\cE$ and $X_{\cF}^+$ are not isomorphic. If the assertion were false, then we could find that $Y_\cF$ and $D_2 (\sigma)$ are isomorphic. The $Y_\cF$ is a hypersurface in $\bP^4$ of degree $5$ by Theorem \ref{mainG24} \eqref{mainG242}. On the other hand, we may regard $D_2 (\sigma)$ as a complete intersection of two hypersurfaces in $\bP^5$ of degree $2$ and $4$. Therefore the Picard group of $Y_\cF$ (resp.\ $D_2 (\sigma)$) is isomorphic to that of $\bP^4$ (resp.\ $\bP^5$) by Lefschetz theorem for Picard groups, see, e.g., \cite[Example 3.1.25]{LazPAGI}. Let the generator of $\Pic (Y_\cF)$ and $\Pic (D_2 (\sigma))$ be $A_1$ and $A_2$ respectively. We deduce the contradiction $5 = A_1^3 = A_2^3 = 8$.
\end{remark}


\subsection{\texorpdfstring{Another description of $\bP (\cF)$}{Another description of P(F)}}

As in Section \ref{V5_sec}, let $W$ be a five-dimensional vector space over $\bC$. Fix a nonzero vector $w_0 \in W$ and let $V=W/\bC w_0$.
Then there is a natural morphism 
\begin{equation}\label{G24cok}
    0 \to \cO_{\bP(W^{\vee})} (-1) \ra V \otimes \cO_{\bP(W^{\vee})}
\end{equation}
via the Euler sequence
\begin{equation*}
    0 \to \cO_{\bP(W^{\vee})} (-1) \to W \otimes \cO_{\bP(W^{\vee})} \to T_{\bP (W^\vee)} (- 1) \to 0
\end{equation*}
and the quotient map $W \twoheadrightarrow V$. We define a coherent sheaf $\cC' $ on $\bP(W^{\vee})$ to be the cokernel of \eqref{G24cok}. Then $\cC'$ is locally free of rank three on $\bP (W^{\vee}) \setminus \{[w_0^{\vee}]\}$.

Recall that $0 \ra \cS \ra V \otimes \cO_{\Gr(2,V)} \ra \cQ \ra 0$ is the universal exact sequence on $\Gr(2,V)$ and define $\widetilde{\cS}$ to be the kernel of the surjection
\begin{equation}\label{G24St}
    W \otimes \cO_{\Gr(2,V)} \twoheadrightarrow V \otimes  \cO_{\Gr(2,V)} \twoheadrightarrow \cQ .
\end{equation}
Hence we have an exact sequence
\begin{equation}\label{eulseqTSGr24}
    0 \ra \widetilde{\cS} \ra W \otimes \cO_{\Gr(2,V)} \ra \cQ \ra 0.
\end{equation}
Note that we have $\widetilde{\cS} \cong \cS \oplus \cO_{\Gr(2,V)}$ by choosing a splitting $W \cong V \oplus \bC$. Because $\cS$ has rank two and $\det (\cS) = \cO_{\Gr(2,V) } (-1)$, we find that $\cS \cong \cS^{\vee} \otimes \det (\cS) = \cS^{\vee} (- 1)$ and thus $\tcS^{\vee} (1) \cong \cF$. In particular, we get $\bP (\tcS^{\vee}) \cong \bP (\cF)$ over $\Gr (2, V)$.

By the surjection $ W^{\vee} \otimes \cO_{\Gr(2,V)} \twoheadrightarrow \widetilde{\cS}^{\vee}$, we have a closed immersion
\begin{equation*}
    \bP_{\Gr(2,V)}( \widetilde{\cS}^{\vee}) \hookrightarrow  \bP_{\Gr(2,V)} ( W^{\vee} \otimes \cO_{\Gr(2,V)} ) = \bP(W^{\vee}) \times \Gr(2,V).
\end{equation*}
Let $p_{\bP(W^{\vee}) }$ and $p_{\Gr(2,V)}$ be the natural projections from $\bP(W^{\vee}) \times \Gr(2,V)$ to the first and the second factors respectively. Set $p_1$ be the restriction of $p_{\bP(W^{\vee}) }$ to $\bP_{\Gr(2,V)}( \widetilde{\cS}^{\vee})$, and similarly for $p_2$. Note that 
\begin{equation}\label{G24phiD}
    \cO_{\tcS^{\vee}} (1) = p_1^{\ast} \cO_{\bP (W^{\vee})} (1).
\end{equation}
The following gives another description of the projective bundle associated to $\widetilde{\cS}^{\vee}$.

\begin{lemma}\label{G24anoGRp}
The Grassmannian $G_{\bP(W^{\vee}) }(\cC',2)$ coincides with the projective bundle $\bP_{\Gr(2,V)}( \widetilde{\cS}^{\vee})$ in $\bP(W^{\vee}) \times \Gr(2,V)$: \begin{equation*}
    \begin{tikzcd}
        & \Gr_{\bP (W^\vee)} (\cC', 2) = \bP_{\Gr(2,V)} ( \tcS^{\vee} ) \ar[rd, "p_2"] \ar[ld, "p_1", swap]  &  \\
        \bP(W^{\vee}) && \Gr(2,V).      
    \end{tikzcd}
\end{equation*}
In particular,
$p_1$ is a $\bP^2$-bundle $\bP(\wedge^2 \cC')$ over $\bP(W^{\vee} )\setminus \{[w_{0}^{\vee}]\}$,
and the fiber $p_1^{-1}([w_0^{\vee}]) $ is $\Gr(V, 2 )=\Gr(2,V)$.
\end{lemma}

\begin{proof}
By \eqref{G24St} and Lemma \ref{eq_as zero locus}, the projective bundle $\bP_{\Gr(2,V)}( \widetilde{\cS}^{\vee})$ is the zero locus of
\begin{align}\label{eq_zero_locus}
p_{\Gr(2,V)}^* \cQ^{\vee} \ra W^{\vee}  \otimes \cO_{\bP(W^{\vee}) \times \Gr(2,V)} \ra p_{\bP(W^{\vee}) }^* \cO_{ \bP(W^{\vee}) } (1)
\end{align}
in $\bP(W^{\vee}) \times \Gr(2,V)$. Since the zero locus of \eqref{eq_zero_locus} coincides with that of 
\[
    p_{\bP(W^{\vee}) }^* \cO_{ \bP(W^{\vee}) } (-1)   \ra V  \otimes \cO_{\bP(W^{\vee}) \times \Gr(2,V)} \ra p_{\Gr(2,V)}^* \cQ,
\]
the bundle $\bP_{\Gr(2,V)}( \widetilde{\cS}^{\vee}) $ coincides with
\[
\Gr_{\bP(W^{\vee}) }(\cC',2) \hookrightarrow \Gr_{\bP(W^{\vee}) }(V  \otimes \cO_{ \bP(W^{\vee}) },2) =  \bP(W^{\vee}) \times \Gr(2,V)
\]
again by \eqref{G24cok} and Lemma \ref{eq_as zero locus}.

To see the last statement,
we note that $\cC'$ is locally free of rank $3$ on $\bP(W^{\vee} )\setminus \{[w_{0}^{\vee}]\}$.
Hence the Pl\"ucker embedding 
\[
\Gr_{\bP (W^\vee)} (\cC', 2)
\hookrightarrow \bP(\wedge^2 \cC')
\]
is an isomorphism over $\bP(W^{\vee} )\setminus \{[w_{0}^{\vee}]\}$.
On the other hand, we have
$p_1^{-1}([w_0^{\vee}]) = \Gr(\cC'([w_0^{\vee}]), 2) = \Gr(V, 2 )$.
\end{proof}


\subsection{Construction of the flop}

Recall that the defining section $s_{\sigma}$ of the Calabi--Yau $X_{\cF} \subseteq \bP (\cF) \cong \bP (\tcS^{\vee})$ is induced by a general morphism $\sigma \colon \cE^{\vee} = \cO^{\oplus 3} \to \cF$. From \eqref{eulseqTSGr24}, we have the surjection
\begin{equation}\label{G24suj}
    W^{\vee} \otimes \wedge^2 V = H^0 (W^{\vee} \otimes \cO_{\Gr (2, V)} (1))  \twoheadrightarrow H^0 (\tcS^{\vee} (1)) = H^0 (\cO_{\tcS^{\vee}} (1) \boxtimes \cO_{\Gr (2, V)} (1)).
\end{equation}
As in the proof of Proposition \ref{V5CY3}, there exists a general three-dimensional subspace 
\begin{equation*}
    U_{\sigma} \subseteq W^{\vee} \otimes \wedge^2 V
\end{equation*}
such that the image of $U_{\sigma}$ under \eqref{G24suj} corresponds to the subspace induced by the section $s_{\sigma}$.
Hence 
we can identify 
\[
s_{\sigma} \colon p_{2}^{\ast} \cE^{\vee}= \cO_{\bP (\tcS^{\vee})}^{\oplus 3}  \to \cO_{\cF} (1) = \cO_{\tcS^{\vee}} (1) \otimes p_{2}^{\ast} \cO_{\Gr(2,V)}(1)
\]
with the composite
\begin{equation} \label{G24XF+0}
U_{\sigma} \otimes \cO_{\tcS^{\vee}} (-1)  \hookrightarrow   U_{\sigma} \otimes W \otimes \cO_{\bP (\tcS^{\vee})}  \to \wedge^2 V \otimes  \cO_{\bP (\tcS^{\vee})} \to p_{2}^{\ast} \cO_{\Gr(2,V)}(1)
\end{equation}
of the natural morphisms up to the twist by $ \cO_{\tcS^{\vee}} (1)$.

On the other hand,
the subspace $U_{\sigma}$ also induces a morphism
\begin{equation*} 
U_{\sigma} \otimes \cO_{\bP(W^{\vee})}(-1) \hookrightarrow   U_{\sigma} \otimes W \otimes \cO_{\bP(W^{\vee})}  \to \wedge^2 V \otimes  \cO_{\bP(W^{\vee})}
\end{equation*}
on $\bP(W^{\vee})$.
By composing this with 
$\wedge^2 V  \otimes \cO_{ \bP(W^{\vee})} \twoheadrightarrow  \wedge^2 \cC'$,
which is induced by the natural surjection $ V  \otimes \cO_{ \bP(W^{\vee})} \twoheadrightarrow \cC'$,
we obtain a morphism
\begin{equation} \label{G24XF+}
    \tau \colon U_{\sigma} \otimes \cO_{\bP(W^{\vee})}(-1) \to \wedge^2 \cC'.
\end{equation}

We set $\bP^{\circ} = \bP(W^{\vee} ) \setminus \{[w_0^{\vee}]\}$.
By Lemma \ref{G24anoGRp},
we have an identification
\begin{equation} \label{G24XF+3}
p_1^{-1} (\bP^{\circ}) =\Gr_{\bP^{\circ}} (\cC'|_{\bP^{\circ}}, 2) = \bP_{\bP^{\circ}} (\wedge^2 \cC'|_{\bP^{\circ}}).
\end{equation}
Now we can give another description of $X_{\cF} \cap p_1^{-1} (\bP^{\circ}) $:

\begin{proposition}\label{G24CY3}
The Calabi--Yau threefold $X_{\cF}  \subseteq \bP_{\Gr (2,4)}(\tcS^{\vee}) = \Gr_{\bP(W^{\vee})} (\cC', 2)$ coincides with 
\[
    \bP_{\bP(W^{\vee})} (\coker \tau) \subseteq  \bP_{\bP(W^{\vee})} (\wedge^2 \cC')
\]
on the open subset $p_1^{-1}(\bP^{\circ}) \subseteq \Gr_{\bP(W^{\vee}) }(\cC',2)$
under the identification \eqref{G24XF+3}.
\end{proposition}

\begin{proof}

Applying Lemma \ref{detzer} (1) to $\tau|_{\bP^{\circ} }$,
we see that $\bP_{\bP(W^{\vee})} (\coker \tau)$ coincides with
the zero locus of 
\begin{equation*}
s_{\tau|_{\bP^{\circ} }} \colon  U_{\sigma}\otimes p_1^* \cO_{\bP^{\circ}} (- 1) \to p_1^* (\wedge^2 \cC|_{\bP^{\circ}}')  \to \cO_{ \wedge^2 \cC|_{\bP^{\circ}}'}(1)
\end{equation*}
on $p_1^{-1} (\bP^{\circ}) =\Gr_{\bP^{\circ}} (\cC'|_{\bP^{\circ}}, 2) = \bP_{\bP^{\circ}} (\wedge^2 \cC'|_{\bP^{\circ}})$.
By the definition of the Pl\"ucker embedding,
the $\cO_{ \wedge^2 \cC|_{\bP^{\circ}}'}(1)$ coincides with $ p_2^* \cO_{\Gr(2,V)}(1)|_{p_1^{-1} (\bP^{\circ})}$.
Furthermore, we have $  p_1^* \cO_{\bP(W^{\vee})}(1)=\cO_{ \tcS^{\vee}}(1)$ by \eqref{G24phiD}.
Then $s_{\tau|_{\bP^{\circ} }} $ is nothing but 
the restriction of \eqref{G24XF+0} on $p_1^{-1} (\bP^{\circ})$.
Since \eqref{G24XF+0} can be identified with $s_{\sigma}$,
the zero locus of $s_{\tau|_{\bP^{\circ} }}$ coincides with 
$X_{\cF} =Z(s_{\sigma})$ in $p_1^{-1} (\bP^{\circ})$
and the proposition follows.
\end{proof}

Using the above description of $X_\cF$, we can construct the other small contraction of it as follows.

\begin{proposition}\label{G24YF}
Let $Y_{\cF} \subseteq \bP(W^{\vee})$ be the image $p_1 (X_\cF)$ and $\varphi \coloneqq p_1|_{X_\cF{}} \colon X_\cF \to Y_\cF$. Then $\varphi$ is a small contraction and $Y_{\cF}$ is a quintic hypersurface in $\bP (W^{\vee})$.
\end{proposition}

\begin{proof}
Recall that the coherent sheaf $\wedge^2 \cC'$ has rank three. For $x \in \bP^{\circ}=\bP(W^{\vee} )\setminus \{[w_{0}^{\vee}]\}$, the fiber $\varphi^{-1} (x)$ is $\bP_{\bP(W^{\vee})} (\coker \tau (x)) \cong \bP^{3 - k - 1} = \bP^{2 - k}$ for $k = \rank \tau (x) = 1, 2$ by Proposition \ref{G24CY3}. 
On the other hand, the fiber $\varphi^{-1} ([w_0^{\vee}])$ is a smooth conic, which is a linear section $p_1^{-1}([w_0^{\vee}]) =\Gr (2,V)$ of codimension three, since $U_{\sigma}$ is general. 
Note that the codimension of $D_1 (\tau|_{\bP^{\circ}}) \subseteq \bP^{\circ}$ is $4$ and $D_0 (\tau|_{\bP^{\circ}})=\varnothing $. Therefore, the $\varphi$ is small.

Note that $Y_\cF \setminus \{[w_0^{\vee}]\}$ is the zero locus of $\det (\tau|_{\bP^{\circ}})$, 
which is a global section of the line bundle $\det (\wedge^2 \cC') \otimes \cO_{\bP (W^{\vee})} (3)|_{\bP^{\circ}} \cong \cO_{\bP (W^{\vee})} (5)|_{\bP^{\circ}}$. Then it is the restriction of a quintic hypersurface on $\bP^{\circ}$.
Since $Y_\cF $ is irreducible, we see that $Y_\cF $ is a quintic hypersurface.
In particular,
the $Y_{\cF} $ is Cohen–Macaulay.
Since the singular locus of $Y_{\cF} $ is zero-dimensional,
$Y_{\cF} $ is normal and hence $\varphi$ is a small contraction.
\end{proof}

In the remainder of this subsection, we are going to find the flop of $\varphi$. Let $\mu$ be the blow-up $\tbP \to \bP (W^{\vee})$ at the point $[w_0^{\vee}]$ and $E$ the exceptional divisor of $\mu$. This $\tbP$ is the graph of the rational map $\bP (W^{\vee}) \dashrightarrow \bP (V^{\vee})$ given by the projection from the point $[w_0^{\vee}]$, and we get the morphism $f \colon \tbP \to \bP (V^{\vee})$. Set $\cT' = T_{\bP (V^{\vee})} (- 1)$.

\begin{lemma}\label{G24TC} There is an isomorphism $f^{\ast} \cT' \cong \cC'$ on $\tbP \setminus E \cong \bP^{\circ}$ 
\end{lemma}
\begin{proof} 
On $\tbP \setminus E \cong \bP^{\circ}$,  there is a diagram 
\begin{equation*}
    \begin{tikzcd}
        \mu^*\cO_{\bP (W^{\vee})} (- 1) \ar[d,"\alpha"] \ar[r] & \cO_{\tbP } \otimes V \ar[r] \ar[d, equal] & \mu^*\cC'=\cC' \ar[r] \ar[d,"\beta"] & 0\\
       f^*\cO_{\bP (V^{\vee})} (- 1) \ar[r]& \cO_{\tbP} \otimes V \ar[r]  & f^*\cT' \ar[r] & 0 
    \end{tikzcd}
\end{equation*}
where the top row is from \eqref{G24cok}, and the bottom row is the pullback of the Euler sequence on $\bP(V^\vee)$. The map $\alpha$ is induced by the natural projection $W\ra V$ and $\beta$ is the induced map from the exact sequences. 
Since $\alpha$ is an isomorphism on $\tbP \setminus E $, 
so is $\beta$ and this lemma follows.
\end{proof}

We write $\cU_{\sigma}$ for the trivial bundle $U_{\sigma} \otimes \cO_{\bP (W^{\vee})}$.
From the observation 
$W^{\vee} \otimes \wedge^2 V = H^0(\mu^* \cO_{\bP(W^{\vee})} (1)) \otimes H^0(f^* (\wedge^2 \cT'))$, 
the subspace $U_{\sigma} \subseteq W^{\vee} \otimes \wedge^2 V$ also induces a morphism
\begin{equation}\label{G24rho}
    \varrho \colon \mu^{\ast} (\cU_{\sigma} (- 1)) \to f^{\ast} (\wedge^2 \cT'),
\end{equation}
which coincides with $\tau \colon \cU_{\sigma} (- 1) \to \wedge^2 \cC'$ on $\tbP \setminus E \cong \bP^{\circ}$
under the isomorphism $f^{\ast} \cT' \cong \cC'$ in Lemma \ref{G24TC}. 
Set $\tY_{\cF} = D_2 (\varrho)$. By Lemma \ref{G24TC}, the morphism $\mu$ maps $\tY_{\cF}$ onto $Y_\cF$. As in Section \ref{Nefcone_subsec}, let
\[
    \tX_\cF \subseteq \bP_{\tbP} (\mu^{\ast} (\cU_{\sigma}^{\vee} (1))) \cong \bP(U_{\sigma}^{\vee}) \times \tbP \cong \bP^2 \times \tbP
\] 
be the zero locus of the global section $s_{\varrho^{\vee}}$ of
\begin{equation}\label{G24rhod}
    \cO_{\mu^* \cU_{\sigma}^{\vee}(1)}(1) \boxtimes f^* (\wedge^2 \cT') \cong \cO_{\bP^2}(1)  \boxtimes (f^* (\wedge^2 \cT') \otimes \mu^*\cO(1))
\end{equation}
induced by the dual morphism $\varrho^{\vee}$. Therefore, we find the determinantal contraction 
$\tX_\cF \to \tY_\cF \subseteq \tbP$, which is small, for a general $U_\sigma$ since $\mu^{\ast} (\cU_{\sigma} ( 1)) \otimes f^{\ast} (\wedge^2 \cT')$ is globally generated.


We write $X_{\cF}^{\dagger} \subseteq \bP^2 \times \bP (W^{\vee})$ for the image $(\id_{\bP^2} \times \mu) (\tX_\cF)$. Via the projections, we have the morphism $\varphi^{\dagger} \colon X_{\cF}^{\dagger} \to Y_\cF$ and following commutative diagram
\begin{equation} \label{G24comdig}
    \begin{tikzcd}
       \bP^2 \times \tbP \cong \bP_{\tbP} (\mu^{\ast} (\cU_{\sigma}^{\vee} (1))) \ar[r, symbol=\supseteq] \ar[d, two heads, "\id_{\bP^2} \times \mu", swap, shift right=3.6em] & \tX_\cF  \ar[r, two heads] \ar[d, two heads, "\id_{\bP^2} \times \mu", swap] & \tY_\cF  \ar[d, two heads, "\mu"] \\
        \bP^2 \times \bP (W^{\vee}) \cong \bP_{\bP (W^\vee)} (\cU_{\sigma}^{\vee} (1)) \ar[r, symbol=\supseteq] & X_{\cF}^{\dagger} \ar[r, two heads, "\varphi^{\dagger}"] & Y_\cF .
    \end{tikzcd}
\end{equation}


\begin{proposition}\label{G24flop}
The morphism $\varphi^{\dagger} \colon X_{\cF}^{\dagger} \to Y_\cF$ is small. In particular, $(\varphi^{\dagger})^{- 1} \circ \varphi$ is an isomorphism in codimension one. 
\end{proposition}

\begin{proof}
The morphism $\varphi^{\dagger}$ is small over $Y_\cF \setminus \{[w_0^\vee]\}$ since so is $\tX_\cF \to \tY_\cF$. 
Hence it suffices to show that the fiber of $ \varphi^{\dagger}$ over $[w_0^\vee]$ is one-dimensional.
By construction, the fiber over the point $[w_0^\vee]$ is $X^{\dagger}_{\cF} \cap (\bP^2 \times [w_0^\vee])$, 
which is the image of $\tX_\cF \cap (\bP^2 \times E)=Z(s_{\varrho^{\vee}}) \cap (\bP^2 \times E)$. 
Under the isomorphism $f|_E \colon E \to \bP (V^\vee)$ and \eqref{G24rhod}, the restriction of the section $s_{\varrho^{\vee}}$ to $\bP^2 \times E$ corresponds to a global section $\omega$ of $\cO_{\bP^2} (1) \boxtimes \wedge^2 \cT'$ on $ \bP^2 \times \bP(V^{\vee})$. Note that we can regard the morphism $\omega$ as an element of
\[
    H^0 (\cO_{\bP^2} (1) \boxtimes \wedge^2 \cT') = H^0 (\wedge^2 V \otimes \cO_{\bP^2} (1)).
\]
Furthermore, $\omega$ is general in $H^0 (\wedge^2 V \otimes \cO_{\bP^2} (1)) = \wedge^2 V \otimes U_{\sigma}^{\vee}$
since $\omega$ corresponds to the composite  $U_{\sigma} \hookrightarrow W^{\vee} \otimes \wedge^2 V  \to (\bC w_0)^{\vee} \otimes \wedge^2 V $
and $U_{\sigma}$ is general.

Let $pr_i$ be the projection of $\bP^2 \times \bP (V^\vee)$ on the $i$th factor and $Z(\omega) \subseteq \bP^2 \times \bP(V^{\vee})$ the zero locus of $\omega \in H^0 (\cO_{\bP^2} (1) \boxtimes \wedge^2 \cT')$. 
Under the identification $ E= \bP(V^{\vee})$ by $f|_E$,
we have $\tX_\cF \cap (\bP^2 \times E) = Z(s_{\varrho^{\vee}}) \cap (\bP^2 \times E) = Z(\omega) $
and $X^{\dagger}_{\cF} \cap (\bP^2 \times [w_0^\vee]) = pr_1(Z(\omega))$.
Hence the rest is to show $\dim pr_1(Z(\omega)) =1$.

We claim that $pr_1 (Z(\omega)) \subseteq \bP^2$ is the zero locus of $\wedge^2 \omega \in H^0 (\wedge^4 V \otimes \cO_{\bP^2} (2))$, the Pfaffian of $\omega \in H^0 (\wedge^2 V \otimes \cO_{\bP^2} (1))$.
In fact,
a point $[u^\vee] \in \bP(U_{\sigma}^{\vee}) \cong \bP^2$ is contained in $pr_1(Z(\omega))$
if and only if there exists $[v^\vee] \in \bP(V^{\vee})$ such that $([u^{\vee}],[v^{\vee}] ) \in Z(\omega)$,
that is,
\[
\bC u \hookrightarrow U_{\sigma} \xrightarrow{\omega} \wedge^2 V \to \wedge^2 \cT'([v^{\vee}]) =\wedge^2 (V/\bC v) 
\]
is the zero map.
This condition is equivalent to say that the two-form $\omega ([u^\vee]) \in \wedge^2 V$ is degenerate,
that is,
$\wedge^2 \omega ([u^\vee]) =0$.

Since $\omega \in H^0 (\wedge^2 V \otimes \cO_{\bP^2} (1))$ is general,
$\wedge^2 \omega \in H^0 (\wedge^4 V \otimes \cO_{\bP^2} (2))$ is non-zero.
Hence
$pr_1 (Z(\omega)) =Z(\wedge^2 \omega) \subseteq \bP^2$ is one-dimensional
and this proposition follows.
\end{proof}











By Proposition \ref{G24flop}, 
we have the following proposition.

\begin{proposition}\label{G24XF+mop}
Let $X_{\cF}^{+}$ be the normalization of $X_{\cF}^{\dagger}$. Then $\varphi^+ \colon X_{\cF}^{+} \to Y_\cF$ is the flop of $\varphi \colon X_\cF \to Y_\cF$ and the natural projection gives rise to $X_{\cF}^+ \to \bP^2$.
\end{proposition}

\begin{proof}
Since $X_{\cF}^{+} \dashrightarrow X_{\cF}$ is small and the Picard number of $X_{\cF}^{+}$ is at least two,
$X_{\cF}^{+}$ is $\bQ$-factorical.
Hence $X_{\cF}^{+} \to Y_{\cF}$ is the flop of $\varphi \colon X_\cF \to Y_\cF$.
\end{proof}

\subsection{Birational models}

We note that the flop $X_{\cF}^+$ has the same Picard number as $X_\cF$. There are two natural divisors on it. One is the pullback of $\cO_{\bP(W^{\vee}) } (1) $ via $X^+_{\cF}  \ra \bP(W^{\vee})$, denoted by $H'$. The other one, denoted by $L'$, is the pullback of $\cO_{\cU_{\sigma}^{\vee} (1)}(1)$ via $X_\cF^+ \to X_\cF^{\dagger} \subseteq \bP_{\bP (W^\vee)} (\cU_{\sigma}^{\vee} (1))$. Set $\theta \coloneqq (\varphi)^{- 1} \circ \varphi^+$.



\begin{lemma}\label{G24XF+mat}
For the birational map $\theta \colon X_\cF^+ \dashrightarrow X_\cF$, 
the matrix representation
of $\theta_{\ast} : N^1(X^+_{\cF}) \to N^1(X_{\cF}) $
with respect to $\{L', H'\}$ and $\{L, H\}$ is given by
\[
    \begin{bmatrix}
      \theta_\ast
    \end{bmatrix}
=
    \begin{bmatrix}
      5 & 1 \\
      - 6 & - 1 
    \end{bmatrix} .
\]
\end{lemma}

\begin{proof}
By \eqref{G24phiD} and $\cF = \tcS^\vee (1)$, we see that $\theta_\ast H' = L - H$. For the other divisor, we set $\theta_\ast L' = a L + b H$. To apply Lemma \ref{ab}, we need to compute the intersection numbers $L' \cdot H'^2$ and $L'^2 \cdot H'$. Note that $\theta$ is not an isomorphism because only $X_\cF^+$ has a fibration over $\bP^2$. 

From the commutative diagram \eqref{G24comdig} and the fact that  $\tX_\cF$, $X_\cF^\dagger$ and $X_\cF^+$ are birational, we can reduce the computation to the birational model $\tX_\cF$, which is induced by the dual morphism of \eqref{G24rho}. We write $\alpha$ and $\xi$ in $A^1 (\tbP)$ for for the pullbacks of the hyperplane classes on $\bP (V^\vee)$ and $\bP (W^\vee)$ respectively. Then, applying Proposition \ref{LHinterODP}, we find that\footnote{The formula holds without the Calabi--Yau condition \eqref{CYcond}, see \cite[Remark 4.6]{SSW18}.}
\begin{equation}\label{interLH'}
    L'^{k} \cdot H'^{3 - k} = \int_{\tbP} c_{k + 1} (f^\ast (\wedge^2 \cT ') - \mu^\ast (\cU_{\sigma} (- 1))) \cdot \alpha^{3 - k}.
\end{equation}
Note that we have the total Chern class $f^\ast c (\wedge^2 \cT') = 1 + 2 \alpha + 2 \alpha^2$ and the Chow ring $A (\tbP) = \bZ [\alpha, \xi] / (\alpha^4, \xi^2 - \alpha \xi)$ by the Euler sequence of $\bP (V^\vee)$ and $\tbP = \bP_{\bP (V^\vee)} (\cO (1) \oplus \cO)$ respectively. It is easily seen that 
$L' \cdot H'^2 = 14$ and $L'^2 \cdot H' = 28$. By Lemma \ref{ab} and Table \ref{intAll}, we get that
\[
    \begin{cases}
      28 = 40 a^2 + 48 a b + 13 b^2,\\
      14 = 16 a + 11 b.\\
    \end{cases}
\]
Then the only solution is $(a, b) = (5, -6)$ because $ab < 0$, which completes the proof.
\end{proof}

We can now prove Theorem \ref{mainG24}.

\begin{proof}[Proof of Theorem \ref{mainG24}]
First, we claim that $X_\cE \to \bP^2$ is an elliptic fibration. To see this, we let $X_\cE \to P \to \bP^2$ be the Stein factorization. By the adjunction formula and generic smoothness, we know that a general fiber of $X_\cE \to \bP^2$ is a disjoint union of smooth elliptic curves and is linearly equivalent to the cycle $L^2_\cE$. Let $d$ be the degree of $P \to \bP^2$ and $C$ a connected component of the general fiber. 
Since $C \subseteq X_{\cE} \subseteq \Gr(2,4) \times \bP^2$,
we get the closed embedding
\[
    C \xrightarrow{\sim} \pi_\cE (C) \subseteq D_2 (\sigma) \subseteq \Gr (2, 4) \hookrightarrow \bP^5.
\]
Notice that $d (C \cdot H_\cE) = L^2_\cE \cdot H_\cE = 5$ as we have seen in Proposition \ref{movXE}. 
On the other hand, $C \cdot H_\cE = \pi_{\cE} (C) \cdot \cO_{\bP^5} (1)\geqslant 3 $ since $\pi_\cE (C) \subseteq \bP^5 $ is an elliptic curve.
Hence $d = 1$ and $P \to \bP^2$ is a birational finite morphism, so it is an isomorphism. Thus $X_\cE \to \bP^2$ has connected fibers. Similarly, $X_\cF^+ \to \bP^2$ is also an elliptic fibration. We note that the natural projection $\bP_{\tbP} (\mu^{\ast} (\cU^\vee_{\sigma} (1))) \to \bP^2$ corresponds to the complete linear system $|\cO_{\mu^{\ast} (\cU^\vee_{\sigma} (1))} (1) \boxtimes \mu^\ast \cO_{\bP (W^\vee)} (- 1)|$ and $L'^3 = 47$, $H'^3 = 5$ by \eqref{interLH'}. Hence a general fiber of $X_\cF^+ \to \bP^2$ is linearly equivalent to the cycle $(L' - H')^2$ and $(L' - H')^2 \cdot H' = 5$. The rest is the same as the case $X_\cE \to \bP^2$.

The \eqref{mainG242} and \eqref{mainG243} now follow from Propositions \ref{G24YF}, \ref{G24XF+mop} and Lemma \ref{G24XF+mat}. The \eqref{mainG241} and \eqref{mainG244} follow from the construction of $\pi_\cF$ and Proposition \ref{movXE}. 
\end{proof}

Now the proof of Theorem \ref{mainthm_intro} is completed
by Theorems \ref{mainV4}, \ref{mainV5} and \ref{mainG24}. 

\bibliographystyle{alpha}

\end{document}